\documentclass[11pt,a4paper]{amsart}
\usepackage{amscd}
\usepackage{amssymb}
\usepackage[centering,text={15cm,19.5cm}]{geometry}
\usepackage{graphicx}
\usepackage{color}
\usepackage[all]{xy}
\usepackage{mathrsfs}
\usepackage{marvosym}

\definecolor{shadecolor}{rgb}{1,0.9,0.7}

\newtheorem{theorem}{Theorem}[section]

\theoremstyle{definition}

\newtheorem{example}[theorem]{Example}

\theoremstyle{remark}
\newtheorem{remark}[theorem]{Remark}

\numberwithin{equation}{section}
\numberwithin{figure}{section}

\newcommand{\bbfamily}{\fontencoding{U}\fontfamily{bbold}\selectfont}
\newcommand{\textbb}[1]{{\bbfamily#1}}

\newcommand {\lfor} {\mbox{\textbb{[}}}
\newcommand {\rfor} {\mbox{\textbb{]}}}

\newcommand{\NN} {\mathbb{N}}
\newcommand{\ZZ} {\mathbb{Z}}
\newcommand{\QQ} {\mathbb{Q}}
\newcommand{\RR} {\mathbb{R}}
\newcommand{\CC} {\mathbb{C}}
\newcommand{\FF} {\mathbb{F}}
\newcommand{\PP} {\mathbb{P}}
\renewcommand{\AA} {\mathbb{A}}
\newcommand{\GG} {\mathbb{G}}

\newcommand {\fob}  {\mathfrak{b}}

\newcommand {\foj}  {\mathfrak{j}}

\newcommand {\fop}  {\mathfrak{p}}
\newcommand {\foq}  {\mathfrak{q}}

\newcommand {\cl}  {\operatorname{cl}}

\newcommand {\conv} {\operatorname{conv}}

\newcommand {\id}  {\operatorname{id}}

\newcommand {\Int}  {\operatorname{Int}}

\newcommand {\lra}  {\longrightarrow}

\renewcommand {\max} {{\operatorname{max}}}

\renewcommand{\O}  {\mathcal{O}}

\renewcommand{\P}  {\mathscr{P}}

\newcommand {\Proj} {\operatorname{Proj}}

\newcommand {\scrS}  {\mathscr{S}}

\newcommand {\sing} {\mathrm{sing}}

\newcommand {\Spec} {\operatorname{Spec}}

\newcommand {\X} {\mathfrak X}

\def\mydate{\ifcase\month \or January\or February\or March\or
April\or May\or June\or July\or August\or September\or October\or 
November\or December\fi \space\number\day,\space\number\year}

\begin{document}

\title[Toric degenerations]{An invitation to toric degenerations}

\author{Mark Gross} \address{UCSD Mathematics,
9500 Gilman Drive, La Jolla, CA 92093-0112, USA}
\email{mgross@math.ucsd.edu}
\thanks{This work was partially supported by NSF grants 0505325 and
0805328.}

\author{Bernd Siebert} \address{Department Mathematik,
Universit\"at Hamburg, Bundesstra\ss e~55, 20146 Hamburg,
Germany}
\email{bernd.siebert@math.uni-hamburg.de}

\maketitle

\tableofcontents

\section*{Introduction.}
In \cite{affinecomplex} we gave a canonical construction of
degenerating families of varieties with effective anticanonical
bundle. The central fibre $X$ of such a degeneration is a union of
toric varieties, glued pairwise torically along toric prime divisors.
In particular, the notion of toric strata makes sense on the central
fiber. A somewhat complementary feature of our degeneration is their
toroidal nature near the $0$-dimensional toric strata of $X$; near
these points the degeneration is locally analytically or in the
\'etale topology given by a monomial on an affine toric variety. Thus
in this local model the central fiber is a reduced toric divisor. A
degeneration with these two properties is called a \emph{toric
degeneration}. The name is probably not well-chosen as it suggests a
global toric nature, which is not the case as we will emphasize below.
A good example to think of is a degeneration of a quartic surface in
$\PP^3$ to the union of the coordinate hyperplanes. More generally,
any Calabi-Yau complete intersection in a toric variety has toric
degenerations \cite{GBB}. Thus the notion of toric degeneration is a
very versatile one, conjecturally giving all deformation classes of
Calabi-Yau varieties with maximally unipotent boundary points.

Our construction has a number of remarkable features. It generalizes
the construction of a polarized toric variety from an integral
polyhedron (momentum polyhedron) to non-toric situations in a highly
non-trivial, but canonical fashion. It works order by order, each step
being controlled by integral affine (``tropical'') geometry. This has
the striking consequence that any complex geometry feature of the
degeneration that is determined on a finite order deformation of the
central fibre can, at least in principle, be read off tropically. An
important ingredient in this algorithm is the scattering construction,
introduced by Kontsevich and Soibelman in a rigid-analytic setup in
dimension two \cite{ks}. 

The purpose of these notes is to provide an extended introduction to
\cite{affinecomplex}. The emphasis is on highlighting some features of
the construction by going through examples explicitly. To avoid
repeating ourselves, we will introduce most concepts in an ad hoc
fashion and refer to \cite{affinecomplex} for the general case and
more technical definitions.

\section{Purely toric constructions}

\subsection{Toric varieties from polyhedra}
\label{subsect.ToricGeom} To start with let us recall the
algebraic-geometric construction of a toric variety from a convex
integral polyhedron $\sigma \subseteq M_\RR$, the intersection of
finitely many closed halfspaces. To keep track of functorial behaviour
we work in $M_\RR:=M\otimes_\ZZ\RR$ for some free abelian group
$M\simeq\ZZ^n$ of rank $n$. If $\sigma$ is bounded it is the convex
hull in $M_\RR$ of finitely many points in $M\subseteq M_\RR$. In any
case, for each face $\tau\subseteq\sigma$ we have the cone generated
by $\sigma$ relative to $\tau$:
\[
K_\tau\sigma:= \RR_{\ge0} \cdot(\sigma-\tau)
=\big\{m\in M_\RR\,\big|\, \exists m_0\in \tau, m_1\in\sigma,
\lambda\in\RR_{\ge0}: m=\lambda\cdot(m_1-m_0)  \big\}.
\]
These cones are finite rational polyhedral, that is, there exist $u_1,\ldots,
u_s\in M$ with
\[
K_\tau\sigma=C(u_1,\ldots, u_s):=
\RR_{\ge0}\cdot u_1+\ldots +\RR_{\ge0}\cdot u_s.
\]
Note also that $K_\tau\sigma\cap (-K_\tau\sigma)= T_\tau =\tau-\tau$,
the tangent space to $\tau$. So $K_\tau\sigma$ is strictly convex
if and only if $\tau$ is a vertex.

The integral points of $K_\tau\sigma$ define a subring
$\CC[K_\tau\sigma \cap M]$ of the Laurent polynomial ring
$\CC[M]\simeq\CC[z_1^{\pm1},\ldots,z_n^{\pm n}]$ by restricting the
exponents to integral points of $K_\tau\sigma$. For $m\in
K_\tau\sigma\cap M$ we write $z^m$ for the corresponding monomial. The
invertible elements of this ring are precisely the monomials $cz^m$
with $m\in \Lambda_\tau:= T_\tau\cap M$ and $c\in\CC\setminus \{0\}$. 

\begin{example}\label{expl.WP2}
Let $\sigma\subseteq\RR^2$ be the triangle with vertices $v_1=(0,0)$,
$v_2=(1,0)$, $v_3=(0,a)$ with $a\in\NN\setminus\{0\}$.
Then
\begin{eqnarray*}
K_{v_1}\sigma&=& C\big((1,0),(0,1) \big)\\
K_{v_2}\sigma&=& C\big((-1,a),(-1,0) \big)\\
K_{v_3}\sigma&=& C\big((0,-1),(1,-a) \big)
\end{eqnarray*}
For $i=1,3$ the monoid of integral points
$K_{v_i}\sigma\cap \ZZ^2\subseteq \ZZ^2$ is freely generated by the
primitive generators $(1,0), (0,1)$ and $(0,-1),(1,-a)$ of the
extremal rays. In other words,
\[
\NN^2\lra K_{v_1}\sigma\cap\ZZ^2,\quad
(\alpha,\beta)\longmapsto \alpha\cdot (1,0)+\beta\cdot(0,1)
= (\alpha,\beta) 
\]
and
\[
\NN^2\lra K_{v_3}\sigma\cap\ZZ^2,\quad
(\alpha,\beta)\longmapsto \alpha\cdot (0,-1)+\beta\cdot(1,-a)
= (\beta,-\alpha-\beta a) 
\]
are isomorphisms of additive monoids. This shows
\[
\CC[K_{v_i}\sigma\cap\ZZ^2] \simeq \CC[x,y], \quad i=1,3,
\]
as abstract rings. For $i=2$ the integral generators $(-1,a), (-1,0)$
of the extremal rays of $K_{v_2}\sigma$ generate a proper sublattice
of $\ZZ^2$ of index $a=\det\left( \begin{smallmatrix} -1 &-1\\
a&0\end{smallmatrix}\right)$. Thus $(-1,0), (-1,a)$ also do not
suffice to generate $K_{v_2}\sigma\cap\ZZ^2$ as a monoid, for $a>1$.
It is not hard to show that a minimal set of generators of
$K_{v_2}\sigma$ rather consists of the $a+1$ elements
\[
(-1,0),(-1,1),\ldots,(-1,a).
\]
A good way to view $\CC[K_{v_2}\sigma\cap\ZZ^2]$ is as the ring of
invariants of $\CC[x,y]$ under the diagonal action of $\ZZ/a$ by
$a$-th roots of unity $\zeta\in\CC$, $\zeta^a=1$:
\[
x\longmapsto \zeta\cdot x,
\quad
y\longmapsto \zeta\cdot y.
\] 
Under this identification $z^{(-1,i)}\in \CC[K_{v_2}\sigma\cap\ZZ^2]$
corresponds to the invariant monomial $x^i y^{a-i}$.

The remaining rings associated to higher dimensional faces of $\sigma$
are
\begin{eqnarray*}
\CC[K_\tau\sigma]\simeq\begin{cases}
\CC[x,y^{\pm 1}],&\dim\tau=1\\
\CC[x^{\pm 1},y^{\pm 1}],&\tau=\sigma.
\end{cases}
\end{eqnarray*}
\vspace{-4ex}

\qed
\end{example}

As the example indicates, rings of the form $\CC[K_\tau\sigma\cap M]$
(\emph{toric rings}) can be difficult to describe in terms of
generators and relations. To obtain examples that can be easily
written down in classical projective algebraic geometry, in this paper
we therefore almost exclusively restrict ourselves to polyhedra
$\sigma$ with $K_v\sigma\cap M\simeq \NN^n$ as a monoid, for any
vertex $v\in\sigma$. If $m_1,\ldots,m_s$ are the generators of the
extremal rays of $K_v\sigma$, a necessary and sufficient condition for
this to be true is $s=n$ and $\det(m_1,\ldots,m_s)=1$.

Now given a convex integral polyhedron $\sigma\subseteq M_\RR$,
with $\dim\sigma=n$ for simplicity, and a face $\tau\subseteq
\sigma$ we obtain the affine toric variety
\[
U_\tau:=\Spec \big(\CC[K_\tau\sigma\cap M]\big).
\]
Since $\CC[K_\tau\sigma\cap M] \subseteq \CC[M]$ any $U_\tau$ contains
the algebraic torus $U_\sigma=\Spec\big(\CC[M]\big)\simeq \GG_m^n$. More
generally, if $\tau\subseteq\tau'$ then $\CC[K_\tau\sigma\cap M]$ is
canonically a subring of $\CC[K_{\tau'}\sigma\cap M]$, and hence we have
an open embedding
\[
U_{\tau'}\lra U_\tau.
\]
These open embeddings are mutually compatible. Hence the
$U_\tau$ glue to a scheme $X_\sigma$ of dimension $\dim\sigma$. In
other words, there are open embeddings $U_\tau\to X_\sigma$ inducing
the morphisms $U_{\tau'}\to U_\tau$ for all $\tau\subseteq \tau'\subseteq
\sigma$. The multiplication action on $U_\sigma=\GG_m^n$ extends to
$X_\sigma$. Hence $X_\sigma$ is a \emph{toric variety}. Note that
according to this definition toric varieties have a distinguished
closed point, the unit of $\GG_m^n$. Moreover, for faces
$\tau'\subseteq \tau\subseteq \sigma$ the ring epimorphism
\[
\CC[K_{\tau'}\sigma\cap M] \lra \CC[K_{\tau'}\tau\cap M],
\quad
z^m\longmapsto \begin{cases} z^m,& m\in K_{\tau'}\tau\\
0,&\text{otherwise},
\end{cases}
\]
induces a closed embedding $\iota_\tau: X_\tau\to X_\sigma$ with image
disjoint from $\GG_m^n\subseteq X_\sigma$ unless $\tau=\sigma$. The
images of the various $X_\tau$ are called \emph{toric strata} of
$X_\sigma$, the image of $U_\sigma= \GG_m^n$ the \emph{big cell}.
If $\tau,\tau'\subseteq \sigma$ are two faces it holds
\[
\iota_\tau(X_\tau)\cap \iota_{\tau'}(X_{\tau'})= \iota_{\tau\cap\tau'}
(X_{\tau\cap\tau'}).
\]
Here we make the convention $X_\emptyset:=\emptyset$. Hence the face
lattice of $\sigma$ readily records the intersection pattern of the
toric strata of $X_\sigma$. In particular, the facets (codimension one
faces) of $\sigma$ are in one-to-one correspondence with the
\emph{toric prime divisors}, the irreducible Weil divisors that are
invariant under the torus action.

\begin{example}
For $\sigma=\conv\{(0,0),(1,0),(0,a)\}$ from Example~\ref{expl.WP2}
we claim that $X_\sigma$ is the weighted projective plane
$\PP(1,a,1)$. Recall that $\PP(1,a,1)$ is the quotient of
$\AA^3\setminus \big\{(0,0,0)\big\}$ by the action of $\GG_m$ that on
closed points is given by
\[
\lambda\cdot(x_0,x_1,x_2)= (\lambda x_0,\lambda^a x_1, \lambda x_2),
\quad \lambda\in\CC^*.
\]
In fact, on $\AA^3\setminus V(x_0)$ the ring of invariants of the
action is $\CC[x,y]$ with $x=x_1/x_0^a$, $y=x_2/x_0$ and similarly on
$\AA^3\setminus V(x_2)$. On the other hand, on $\AA^3\setminus V(x_1)$
the ring of invariants is generated by $x_0^i x_2^{a-i}/x_1$,
$i=1,\ldots,a$. This is the ring of invariants of the diagonal
$\ZZ/a$-action on $\AA^2$. Hence $\PP(1,a,1)$ has an affine open
covering with spectra of the rings $\CC[K_{v_i}\sigma\cap \ZZ^2]$
discussed in Example~\ref{expl.WP2}. The gluing morphisms between
these open sets are the same as given by toric geometry.
\end{example}

By construction the scheme $X_\sigma$ depends only on the cones
$K_\tau\sigma$, hence only on the \emph{normal fan} of $\sigma$ with
elements the dual cones $(K_\tau\sigma)^\vee\subseteq M_\RR^*$. More
generally, toric varieties are constructed from fans. In particular,
integrality and boundedness of $\sigma$ can be weakened to rationality
of the cones $K_\tau\sigma$. Those toric varieties coming
from integral polyhedra (bounded or not) are endowed with a toric
ample line bundle. In fact, defining the \emph{cone over $\sigma$}
\[
C(\sigma):= \cl\big(\RR_{\ge0}\cdot (\sigma\times\{1\})\big)
\subseteq M_\RR\times\RR,
\]
the ring $\CC[C(\sigma)\cap (M\times\ZZ)]$ is graded by $\deg
z^{(m,h)}:=h \in\NN$. Taking the closure $\cl$ here is important in
the unbounded case. It adds the \emph{asymptotic cone} $\lim_{a\to
0} a\cdot \sigma$ to $M_\RR\times\{0\}$. It is then not hard to see that
one has a canonical isomorphism
\[
X_\sigma\simeq \Proj\big( \CC[C(\sigma)\cap
(M\times\ZZ)] \big).
\]
Although $\CC[C(\sigma)\cap (M\times\ZZ)]$ is not in general generated
in degree $1$, integrality of the vertices of $\sigma$ implies that
the sheaf $\O(1)$ on the right-hand side is nevertheless locally free.
This yields the toric ample line bundle mentioned above.

\subsection{Toric degenerations of toric varieties}
\label{subsect.DegenToric}
Now let us see how certain unbounded polyhedra naturally lead to toric
degenerations with general fibre a toric variety. Let $\tilde\sigma
\subseteq M_\RR\times\RR$ be an $(n+1)$-dimensional convex integral
polyhedron that is closed under positive translations in the last
coordinate: 
\[
\tilde\sigma= \tilde\sigma+ \big(0\times\RR_{\ge0}\big).
\]
Let $q:M_\RR\times\RR\to M_\RR$ be the projection and
\[
\sigma:=q(\tilde\sigma).
\]
Then the non-vertical part of $\partial\tilde\sigma$ is the graph of a
piecewise affine function
\[
\varphi:\sigma\to\RR
\]
with rational slopes. The domains of affine linearity of $\varphi$
define a decomposition $\P$ of $\sigma$ into convex polyhedra. In
terms of this data, $\tilde\sigma$ is the upper convex hull of the
graph of $\varphi$:
\[
\tilde\sigma=\big\{ (m,h)\in M_\RR\times\RR\,\big|\,
h\ge \varphi(m)\big\}.
\]
Thus $\tilde\sigma$ is equivalent to a polyhedral decomposition $\P$
of the convex integral polyhedron $\sigma$ together with a function
$\varphi$ on $\sigma$ that is piecewise affine and strictly convex with
respect to $\P$ and takes integral values at the vertices of $\P$.

Now $X_{\tilde\sigma}$ is an $(n+1)$-dimensional toric variety that
comes with a toric morphism
\[
\pi:X_{\tilde\sigma}\to \AA^1.
\]
In fact, each of the rings $\CC[ K_{\tilde\tau} \tilde\sigma
\cap(M\times\ZZ)]$ is naturally a $\CC[t]$-algebra by letting
$t=z^{(0,1)}$, and the gluing morphisms are homomorphisms of
$\CC[t]$-algebras. The preimage of the closed point $0\in\AA^1$ is
set-theoretically the union of toric prime divisors of
$X_{\tilde\sigma}$ that map non-dominantly to $\AA^1$. It is reduced
if and only if $\varphi$ has integral slopes, that is, takes integral
values at all integral points, not just the vertices. To see this let
$\tilde v=(v,\varphi(v)) \in\tilde \sigma$ be a vertex and
$\varphi_v(m)$ a piecewise linear function on $M_\RR$ which agrees
with $\varphi(v+m)-\varphi(v)$ close to $0$. In other words, the graph
of $\varphi_v(m)$ is the boundary of the tangent cone
\[
K_{\tilde
v}\tilde\sigma= \big\{ (m,h)\in M_\RR\times\RR\,\big|\,
h\ge\varphi_v(m)\big\}
\]
of $\tilde\sigma$ at $\tilde v$. A $\CC$-basis for $\CC[K_{\tilde
v}\tilde\sigma\cap (M\times\ZZ)]/(t)$ is given by
\[
z^{(m,h)},\quad \varphi_v(m)\le h<\varphi_v(m)+1,
\]
and $z^{(m,h)}$ is nilpotent modulo $(t)$ if and only if
$\varphi_v(m)< h$, that is, if $\varphi_v(m)$ is not integral.

\emph{Assume now that $\varphi(m)\in\ZZ$ for all $m\in\sigma\cap
M$.}\, Then $\CC[K_{\tilde v}\tilde\sigma\cap (M\times\ZZ)]/(t)$ has
one monomial generator $z^{(m,\varphi_v(m))}$ for any
$m\in K_v\sigma\cap M$, and the relations are
\[
z^{(m,\varphi_v(m))}\cdot z^{(m',\varphi_v(m'))}=
\begin{cases}
z^{(m+m',\varphi_v(m+m'))},&\exists\tau\in\P:v\in\tau
\text{ and } m,m'\in K_v\tau\\
0,&\text{otherwise}.
\end{cases}
\]
In other words, $\pi^{-1}(0)$ is the scheme-theoretic sum (fibred
coproduct) of the $n$-dimensional toric varieties
\[
X_{\tilde \tau}\simeq X_{q(\tilde \tau)}
\]
with $\tilde\tau\subset \partial\tilde\sigma$ projecting bijectively onto
some $\tau\in\P^{[n]}$. These are precisely the toric prime divisors of
$X_{\tilde\sigma}$ mapping non-dominantly to $\AA^1$. As in
\cite{affinecomplex} $\P^{[k]}$ denotes the set of
$k$-dimensional cells of the polyhedral complex $\P$.

To understand general fibres $\pi^{-1}(t)$, $t\neq 0$, we localize at
$t$. This has the effect of removing the lower boundary of
$\tilde\sigma$, that is, of going over to
$\tilde\sigma+(0\times\RR)=\sigma\times\RR$. Thus
\[
\pi^{-1}(\AA^1\setminus\{0\})= X_{\sigma\times\RR}=X_\sigma\times
(\AA^1\setminus\{0\}).
\]
Thus each general fibre $\pi^{-1}(t)$, $t\neq0$, is canonically
isomorphic to $X_\sigma$. Note, however, that these isomorphisms
degenerate as $t$ approaches $0$.

\begin{example}\label{expl.DegenP1}
Here is a degeneration of $\PP^1$ to two copies of $\PP^1$ featuring
an $A_{l-1}$-singularity in the total space. Let $\sigma=[0,a+1]$,
$\P^{[1]}=\big\{ [0,a],[a,a+1]\big\}$, $\varphi(0)=\varphi(a)=0$,
$\varphi(a+1)=l$ as in Figure~\ref{fig.DegenP1}. The slopes of
$\varphi$ are $0$ and $l$ on the two $1$-cells.

The boundary of $\tilde\sigma$ has two non-vertical components.
Each gives one of the two irreducible components of $\pi^{-1}(0)$.
Their point of intersection is the $0$-dimensional toric stratum
defined by the vertex $\tilde v=(a,0)$ of $\tilde\sigma$. The monoid
$K_{\tilde v}\tilde\sigma\cap \ZZ^2$ has generators $m_1=(-1,0)$,
$m_2=(1,l)$, $m_3=(0,1)$ with relation $m_1+m_2=l\cdot m_3$. Hence,
\[
\CC[K_{\tilde v}\tilde\sigma\cap\ZZ^2]\simeq
\CC[z_1,z_2,t]/(z_1z_2-t^l),
\]
with $z_i=z^{m_i}$ for $i=1,2$ and $t=z^{(0,1)}$ defining $\pi$. This
is a local model of a smoothing of a nodal singularity with an
$A_{l-1}$-singularity in the total space. Thus the changes of slope of
$\varphi$ at the non-maximal cells of $\P$ determine the
singularities of the total space.\\[2ex]
\begin{figure}[ht]
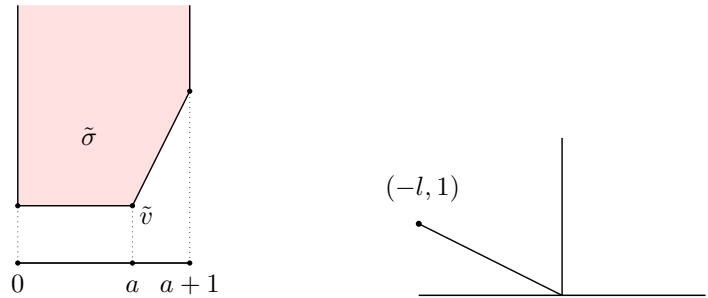

\hfill\input{DegenP1.pstex_t}
\hfill\input{FanDegenP1.pstex_t}
\hfill\mbox{}
\caption{The polyhedron $\tilde\sigma$ (left) and the normal fan of
$\tilde\sigma$ (right).}
\label{fig.DegenP1}
\end{figure}

Another way to understand the total space is from the normal fan of
$\tilde\sigma$. It can be obtained by subdividing the fan of
$\AA^1\times\PP^1$ by the ray through $(-l,1)$. This corresponds to a
weighted blow up at one of the zero-dimensional toric strata, leading
to another $\PP^1$ over $0\in\AA^1$ and the $A_{l-1}$-singularity.

Note also that the length $a$ of the interval is completely irrelevant
to the complex geometry; it only changes the polarization which we
did not care about at this point.
\qed
\end{example}

\begin{example}
\label{expl.degenWP1xP1}
For a two-dimensional example consider the convex hull
\[
\sigma=\conv \big\{(1,0),(0,1),
(-1,0), (0,-1)\big\}
\]
with the polyhedral decomposition into $4$~standard $2$-simplices
shown in the figure below. At $(x,y)\in\sigma$ the function $\varphi$
is given by $0$, $x$, $y$ and $x+y$, respectively, depending on the
maximal cell containing $(x,y)$ as shown in the figure. Thus the slope
of $\varphi$ changes by $1$ along each interior $1$-cell of $\P$.

The associated degeneration has as central fibre a union of $4$~copies
of $\PP^2$, glued pairwise along the coordinate lines as indicated by
$\P$. Thus the singular locus is a union of $4$ projective lines,
joined in one point $X_{\tilde v}$, where $\tilde v=
(0,0,0)\in\tilde\sigma$. The monoid $K_{\tilde v}\tilde\sigma$ is
generated by $m_1=(1,0,1)$, $m_2=(0,1,1)$, $m_3=(-1,0,0)$,
$m_4=(0,-1,0)$ fulfilling the single relation $m_1+m_3=m_2+m_4$. Thus
\[
\CC[K_{\tilde v}\tilde\sigma\cap\ZZ^3]\simeq
\CC[z_1,z_2,z_3,z_4]/(z_1z_3-z_2z_4),
\]
with $\CC[t]$-algebra structure defined by $t=z^{(0,0,1)}
=z_1z_3=z_2z_4$. This shows that the total space $X_{\tilde\sigma}$
has a singular point isomorphic to the origin in the affine cone
over a smooth quadric, while the central fibre is a product of two
normal crossing singularities.
\medskip

\noindent
\begin{minipage}{11cm}
\ \ The general fibre is isomorphic to $X_\sigma$, as always, which here
is a toric $\ZZ/2$-quotient of $\PP^1\times\PP^1$. In fact,
$\left(\begin{smallmatrix} 1&1\\-1&1\end{smallmatrix}\right)$ maps the fan of
$\PP^1\times\PP^1$ to the normal fan of $\sigma$. Restricted to
the big cell $\GG_m^2\subset \PP^1\times\PP^1$, this map is given by
\[
\CC[u^{\pm1},v^{\pm1}]\lra \CC[x^{\pm1},y^{\pm1}],\quad
u\longmapsto xy,\ v\longmapsto x^{-1}y.
\]
The subring of $\CC[x^{\pm 1},y^{\pm 1}]$ generated by $x^{\pm1}
y^{\pm 1}$ is the invariant ring for the involution $(x,y)\to(-x,-y)$,
and this involution extends to $\PP^1\times\PP^1$. Note that
$X_\sigma$ has $4$~isolated quotient singularities. These correspond
to the vertices of $\sigma$.
\end{minipage}
\hfill
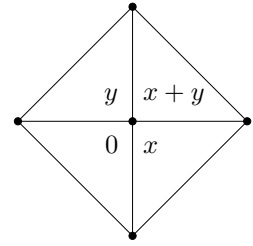
\begin{minipage}{30mm}
\input{DegenWP1xP1.pstex_t}
\end{minipage}
\bigskip

It is also instructive to write down this degeneration embedded
projectively. As explained at the end of \S\ref{subsect.ToricGeom} we
have to take the integral points of the cone $C(\tilde\sigma)
\subseteq \RR^4$ over $\tilde\sigma$ as generators of a graded
$\CC[t]$-algebra, with the degree given by the projection to the last
coordinate and $t=z^{(0,0,1,0)}$. Now $\CC[ C(\tilde\sigma)\cap
\ZZ^4]$ is generated as a $\CC[t]$-algebra by the monomials
\[
X=z^{(1,0,1,1)},\ Y=z^{(-1,0,0,1)},\ Z=z^{(0,1,1,1)},\ 
W=z^{(0,-1,0,1)},\ U=z^{(0,0,0,1)}.
\]
Note that these generators are in one-to-one correspondence
with the integral points of $\sigma$. The relations are
\[
XY-t U^2,\ ZW-t U^2.
\]
This exhibits $X_{\tilde\sigma}$ as the intersection of two quadrics in
$\PP^4_{\AA^1}=\AA^1\times\PP^4$.
\qed
\end{example}

\section{Introducing singular affine structures}
\label{sect.SingAffine}
From a birational classification point of view toric varieties are
boring as they are all rational. In particular, it is impossible to
construct degenerations of non-rational varieties directly by the
method of \S\ref{subsect.DegenToric}. The idea in \cite{affinecomplex}
is that one can get a much larger and more interesting class of
degenerations by gluing toric pieces in a non-toric fashion. The
central fibre is still represented by a cell complex $\P$ of integral
polyhedra, but the integral affine embedding of the cell complex into
$\RR^n$ exists only locally near each vertex. In other words, the
underlying topological space of $\P$ is an integral affine manifold
$B$, with singularities on a cell complex $\Delta\subseteq B$ of real
codimension~$2$ that is a retract of $|\P^{[n-1]}|\setminus
|\P^{[0]}|$. We then construct the total space $\X$ of the
degeneration order by order, by gluing torically constructed
non-reduced varieties, ``thickenings'' of toric varieties so to speak,
in a non-toric fashion.

Thus our starting data are integral cell-complexes with compatible
integral affine charts near the vertices. We call these \emph{integral
tropical manifolds} (\cite{affinecomplex}, Definition~1.2) because
they arise naturally as the bounded parts of the embedded tropical
varieties associated to the degeneration.

\subsection{Degenerations of hypersurfaces}
\label{subsect.DegenHypersurface}
A hypersurface $X\subseteq \PP^{n+1}$ of degree $d\le n+2$ can be
degenerated to a union of $d$ coordinate hyperplanes simply by
deforming the defining equation. For example, for $n=2$ let
$f\in\CC[X_0,\ldots,X_3]$ be a general homogeneous polynomial of
degree $d\le 4$. Then
\[
X_0\ldots X_{d-1}+tf=0
\]
defines a family $\pi: Y\to\AA^1$ with $\pi^{-1}(0)$ a union of $d$
coordinate hyperplanes. This is not a semistable family because $Y$ is
not smooth at the intersection of $V(f)$ with the singular locus of
$V(X_0\cdot\ldots\cdot X_{d-1})$. The latter consists of $\binom{d}{2}$
projective lines $V(X_i,X_j)$, $0\le i<j\le d-1$. Since $f$ is general
the intersection of $V(f)$ with any of these projective lines consists
of $d$ reduced points with two nonzero coordinate entries each, that
is, not equal to the $4$~special points $[1,0,0,0]$, $[0,1,0,0]$,
$[0,0,1,0]$, $[0,0,0,1]$. Near any of these points $Y\to\AA^1$ is
locally analytically given by the projection of the three-dimensional
$A_1$-singularity
\[
V(xy-wt)\subseteq\AA^4
\]
to the $t$-coordinate. Note that for $w\neq0$ this is a product of a
semistable degeneration of a curve with $\AA^1\setminus\{0\}$, but
this fails at $w=0$, which contains the singular point of $Y$.

While the local model of this degeneration is still toric, the
singular points of $Y$ are general points of one-dimensional toric
strata of $\pi^{-1}(0)$. Hence this is a very different degeneration
than the torically constructed ones in \S\ref{subsect.DegenToric}. Our
first aim is to obtain this local degeneration naturally from a
tropical manifold. This was the starting point for \cite{affinecomplex}
in March~2004.

\subsection{A singular affine manifold}
\begin{figure}[ht]
\input{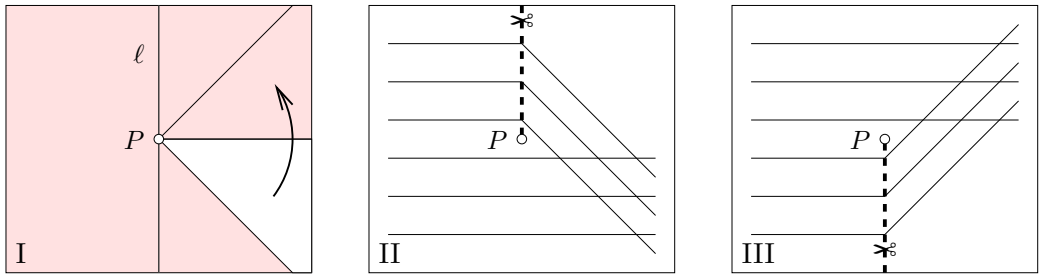}
\caption{Various charts I--III for the focus-focus singularity: On the
left, the shaded area is a fundamental domain for $B$. In the two
other figures, the dashed vertical lines indicate the parts of $\ell$
to be removed; the horizontal curves are all straight affine lines.}
\label{fig.focus-focus}
\end{figure}

There is a famous two-dimensional singular integral affine manifold in
the theory of integrable systems, called the \emph{focus-focus
singularity} \cite{williamson}, which is the model for our singular
affine manifolds $B$ at general points of $\Delta$, the
codimension~$2$ singular locus of the affine structure. An (integral)
affine structure on a topological manifold is an atlas with transition
functions in the (integral) affine linear group. Parallel transport of
tangent vectors is well-defined on such manifolds. In fact, an affine
manifold comes naturally with a flat, torsion-free, but usually
non-metric connection. The focus-focus singularity is (the germ at the
origin $P$ of) $\RR^2$ with an integral affine structure away from $P$
such that parallel transport counterclockwise around $P$ gives the
transformation $\big(\begin{smallmatrix}1&0\\
1&1\end{smallmatrix}\big)$. This affine manifold $B$ can be
constructed by gluing $\RR^2\setminus (\RR_{\ge0}\times\{0\})$ with
the given affine structure to itself by the integral affine
transformation
\[
(x,y)\longmapsto (x,x+y),\quad\text{for}\ x\ge 0,\ -x<y<0,
\]
see Chart~I in Figure~\ref{fig.focus-focus}. Note that $\partial_y$ is
an invariant tangent vector, and indeed the projection $(x,y)\to x$ on
$\RR^2\setminus (\RR_{\ge0}\times\{0\})$ to the first coordinate
descends to a continuous map $B\to \RR$ that is integral affine away
from $P\in B$. The preimage of the origin is a line $\ell$ through
$P$.

Alternatively, $B$ can be constructed from the two charts
$\RR^2\setminus(\{0\} \times\RR_{\ge 0})$ and $\RR^2\setminus(\{0\}
\times\RR_{\le 0})$, each covering $B$ minus one half of~$\ell$, via
\[
\big(\RR\setminus\{0\}\big) \times\RR
\lra \big(\RR\setminus\{0\}\big)\times\RR,
\quad
(x,y)\longmapsto \begin{cases} (x,y),& x<0\\
(x,x+y),&x>0,\end{cases} 
\]
see Charts~II and III in Figure~\ref{fig.focus-focus}. Note that we do
not know how to continue any non-vertical affine line across $P$.

\subsection{The basic example}
\label{subsect.BasicExample}
The focus-focus singularity admits a polyhedral decomposition by
decomposing along the invariant line $\ell$. This decomposition has
two maximal cells $\sigma_\pm$, with preimages $\RR_{\ge 0}\times\RR$
and $\RR_{\le 0}\times\RR$ in any of the charts of
Figure~\ref{fig.focus-focus}. We can also define a strictly convex,
integral affine function $\varphi:B\to\RR$ by letting $\varphi=0$ on
the left maximal cell $\sigma_-$ and $\varphi(x,y)= x$ on the right
maximal cell $\sigma_+$. Again, this is independent of the chosen
chart, and it takes integral values at any integral point of $B$.

Now let us try to treat this situation with a singularity in the same way
as we did in \S\ref{subsect.DegenToric}. We can still describe the central
fibre as before:
\begin{eqnarray*}
&&X_0= X_{\sigma_-}\amalg_{X_\rho} X_{\sigma_+},\\
&&X_{\sigma_-}=\Spec \big(\CC[x, w^{\pm1}]\big),\ 
X_{\sigma_+}=\Spec \big(\CC[y, w^{\pm1}]\big),\ 
X_\rho=\Spec \big(\CC[w^{\pm1}]\big),
\end{eqnarray*}
where $\rho=\sigma_-\cap \sigma_+$, $x=z^{(-1,0)}$, $y=z^{(1,0)}$,
$w=z^{(0,1)}$. By writing this it is important to fix one affine
realization of $\sigma_{\pm}$, say the Chart~II in
Figure~\ref{fig.focus-focus}. Intrinsically, the exponent $(-1,0)$
defining $x$ is an integral tangent vector on $\sigma_-$, while the
exponent $(1,0)$ defining $y$ is an integral tangent vector on
$\sigma_+$. For $w$ we can take either maximal cell, because the
exponent $(0,1)$ just represents the global integral vector field
$\partial_y$ defined above.

Continuing to work in Chart~II the procedure of
\S\ref{subsect.DegenToric} yields the deformation $V(xy-t)$ of
$V(xy)\subseteq \AA^2\times(\AA^1\setminus\{0\})$. But in Chart~III
the tangent vector $(1,0)$ of $\sigma_+$ changes to $(1,1)$. Since
$z^{(1,1)}= wy$ the same procedure applied in this chart therefore
yields $V(xy-wt)$ for the deformation. Because $w$ is invertible
$V(xy-t)$ and $V(xy-wt)$ are isomorphic as schemes over $\CC[t]$, but
not as deformation with the given embedding of $X_0$. Of course, this
phenomenon is just due to the affine monodromy around $P$. As
expected, the singularities of the affine structure lead to
inconsistencies in the naive approach using toric geometry.

The starting point to overcome this problem is to work modulo
$t^{k+1}$. In other words, we want to construct a $k$-th order
deformation $X_k\to \Spec \big(\CC[t]/(t^{k+1})\big)$. The topological
spaces of $X_0$ and $X_k$ are the same, so we only have to deal with
the structure sheaf. We now use toric geometry merely to define the
correct non-reduced versions ($k$-th order thickenings) of the toric
strata. We then glue two maximal-dimensional strata intersecting in
codimension one by a non-toric automorphism of the common codimension
one stratum, but the choice of this automorphism is different
depending on which chart we use.

In the present example, the thickenings of the toric strata
suggested by toric geometry are given by the rings
\begin{eqnarray*}
R^k_{\sigma_-}&=& S_k[x_1,y_1,w^{\pm 1}]/\big(x_1y_1-t,y_1^{k+1}\big),\\
R^k_{\sigma_+}&=& S_k[x_2,y_2,w^{\pm 1}]/\big(x_2y_2-t,x_2^{k+1}\big),\\
R^k_{\rho,\sigma_-}&=&
S_k[x_1,y_1,w^{\pm 1}]/\big(x_1y_1-t,x_1^{k+1},y_1^{k+1}\big),\\
R^k_{\rho,\sigma_+}&=&
S_k[x_2,y_2,w^{\pm 1}]/\big(x_2y_2-t,x_2^{k+1},y_2^{k+1}\big),
\end{eqnarray*}
with $S_k=\CC[t]/(t^{k+1})$. Here we distinguish between $x$ and $y$
as monomials on $\sigma_-$ and on $\sigma_+$. Thus $y_1=z^{(1,0,1)}$ in
either chart, while $y_2=z^{(1,0,1)}$ in Chart~II and $y_2=z^{(1,1,1)}$ in
Chart~III. These monomials only depend on the affine structure on
$\sigma_\pm$ and hence have an intrinsic meaning. For the thickening
of the $\rho$-stratum, however, we obtain two rings, depending on
which maximal cell $\sigma_{\pm}$ the monomials live on. For the
monomial $w$ we don't need to make any choices because it corresponds
to a globally defined vector field. In any case, we have two natural
$S_k$-algebra epimorphisms
\[
q_-:R^k_{\sigma_-}\lra R^k_{\rho,\sigma_-},\quad
q_+:R^k_{\sigma_+}\lra R^k_{\rho,\sigma_+},
\]
exhibiting the thickened $\rho$-stratum as closed subscheme of the
$\sigma_-$- and $\sigma_+$-stratum, respectively.

Affine geometry suggests two isomorphisms $h_{\rm II}, h_{\rm
III}:R^k_{\rho,\sigma_-}\to R^k_{\rho,\sigma_+}$, depending on which
chart we use to go from $\sigma_-$ to $\sigma_+$:
\begin{align*}
\hspace*{20ex} h_{\rm II}:\ &x_1\longmapsto x_2,&&y_1\longmapsto y_2,
&&w\longmapsto w\hspace*{20ex}\\
h_{\rm III}:\ &x_1\longmapsto wx_2,&&y_1\longmapsto w^{-1}y_2,
&&w\longmapsto w
\end{align*}
Now comes the point: To remedy the inconsistency caused by this, let
$\alpha\in\CC^*$ and compose $h_{\rm II}$ with the automorphism
\[
g_{\rm II}:x_2\longmapsto (\alpha+ w)\cdot x_2,\quad
y_2\longmapsto (\alpha+w)^{-1} y_2,\quad
w\longmapsto w
\]
of the localization $\big(R^k_{\rho,\sigma_+}\big)_{\alpha+w}$. For
Chart~III we compose with
\[
g_{\rm III}:x_2\longmapsto (1+\alpha w^{-1})\cdot x_2,\quad
y_2\longmapsto (1+\alpha w^{-1})^{-1} y_2,\quad
w\longmapsto w.
\]
We then indeed obtain
\begin{eqnarray*}
\big(g_{\rm II}\circ h_{\rm II}\big)(x_1)&=& (\alpha+w)\cdot x_2
\ =\ (1+\alpha w^{-1})\cdot w x_2\ =\ 
\big(g_{\rm III}\circ h_{\rm III}\big)(x_1),\\
\big(g_{\rm II}\circ h_{\rm II}\big)(y_1)&=& (\alpha+w)^{-1}\cdot y_2
\ =\ (1+\alpha w^{-1})^{-1}\cdot w^{-1} y_2\ =\
\big(g_{\rm III}\circ h_{\rm III}\big)(y_1).
\end{eqnarray*}
The result of gluing $\Spec R^k_{\sigma_-}$ and $\Spec R^k_{\sigma_+}$
along the codimension one strata via this isomorphism is given by the
fibre product $R^k_{\sigma_-} \times_{(R^k_{\rho,\sigma_+}
)_{\alpha+w}} R^k_{\sigma_+}$. In this fibre product the homomorphism
$R^k_{\sigma_+}\to (R^k_{\rho,\sigma_+} )_{\alpha+w}$ is the
composition of $q_+$ with localization, while $R^k_{\sigma_-}\to
(R^k_{\rho,\sigma_+} )_{\alpha+w}$ composes $q_-$ and the localization
homomorphism with $g_{\rm II}\circ h_{\rm II}= g_{\rm III}\circ h_{\rm
III}$. It can be shown (\cite{affinecomplex}, Lemma~2.34) that
generators for this fibre product as an $S_k[w^{\pm 1}]$-algebra are
\[
x=\big(x_1,(\alpha+w)x_2\big),\
y=\big((\alpha+w)y_1,y_2\big)
\]
with the single relation (coming from $x_1y_1=t=x_2y_2$) 
\begin{equation}\label{eqn.BasicExample}
xy-(\alpha+w)t=0.
\end{equation}
This is the $k$-th order neighbourhood of the surface degeneration
discussed in \S\ref{subsect.DegenHypersurface}, with the
$A_1$-singularity at $x=y=0$ and $w=-\alpha$.

\subsection{General treatment in codimension one}
\label{subsect.DegenCodim1}
The discussion in \S\ref{subsect.BasicExample} generalizes to
arbitrary dimension $n$ as follows. Let $\rho$ be a codimension one
cell of $\P$, separating the maximal cells $\sigma_-$, $\sigma_+$.
This gives three toric varieties $X_{\sigma_\pm}$ and $X_\rho$ with
the respective big cells $\Spec\big( \CC[\Lambda_{\sigma_\pm}]
\big)\simeq \GG_m^n$ and $\Spec\big( \CC[\Lambda_\rho] \big)\simeq
\GG_m^{n-1}$. For an integral polyhedron $\tau$ we use the notation
$\Lambda_\tau$ for the lattice of integral tangent vector fields on
$\tau$. The convex integral piecewise affine function $\varphi$
changes slope along $\rho$ by some integer $k>0$. We then obtain the
same rings and gluing morphisms as in \S\ref{subsect.BasicExample},
except that $w$ has to be replaced by $n-1$ coordinates
$w_1,\ldots,w_{n-1}$ for $\GG_m^{n-1}\subseteq X_\rho$.

There are affine charts on the interiors of $\sigma_\pm$ and near each
vertex of $\rho$. Each vertex $v\in \rho$ suggests an isomorphism
$h_v: R^k_{\rho,\sigma_-} \to R^k_{\rho,\sigma_+}$. For a different
vertex $v'\in \rho$ we have $h_{v'}(z^m)=h_v (z^m)$ for
$m\in\Lambda_\rho$ because these tangent vectors are globally defined
on $\sigma_-\cup\sigma_+$. But for general $m\in\Lambda_{\sigma_-}$,
parallel transport $\Pi_{v'}$ to $\sigma_+$ through $v'$ is related to
parallel transport $\Pi_v$ through $v$ as follows:
\begin{eqnarray}\label{eqn.ParallelTransport}
\Pi_{v'}(m)= \Pi_v (m)
+\langle m, \check d_\rho\rangle \cdot m^\rho_{v'v}.
\end{eqnarray}
Here $m^\rho_{v'v}$ is some element of $\Lambda_\rho$, viewed
canonically as a subset of $\Lambda_{\sigma_+}$, and $\check d_\rho$
is the primitive generator of $\big( \Lambda_{v'}/\Lambda_\rho\big)^*
\simeq \ZZ$ that evaluates positively on tangent vectors pointing from
$\rho$ into $\sigma_+$. This monodromy formula implies
\[
h_{v'}(z^m)= \big( z^{m^\rho_{v'v}}\big)^{\langle m, \check
d_\rho\rangle}\cdot h_v(z^m).
\]
Taking $x=z^m$, $y=z^{-m}$ in the chart at $v$ with $\langle m, \check
d_\rho\rangle=1$ yields $h_v'(x)=z^{m^\rho_{v'v}} \cdot h_v(x)$,
$h_v'(y)=z^{-m^\rho_{v'v}} \cdot h_v(y)$. Thus to define the gluing
invariantly we can compose $h_v$ with an automorphism of
$R^k_{\rho,\sigma_+}$ of the form
\[
g_v: z^m\longmapsto f_{\rho,v}^{-\langle \check d_\rho, m\rangle}
\cdot z^m
\]
for functions $f_{\rho,v}\in \CC[\Lambda_\rho]$, indexed by vertices
$v$ of $\rho$ and fulfilling the change
of vertex formula
\begin{equation}
\label{eqn.ChangeOfVertex}
f_{\rho,v'}= z^{m^\rho_{v'v}}\cdot f_{\rho,v}.
\end{equation}
In the example in \S\ref{subsect.BasicExample} we can think of $v$ as
lying on the upper half of $\ell$, $v'$ on the lower half and
$f_{\rho,v} =1+\alpha w^{-1}$, $f_{\rho,v'}= \alpha+w$. Then
$m^\rho_{v'v}= (0,1)$, $z^{m^\rho_{v'v}}= w$, and the change of vertex
formula reads $\alpha+w= w\cdot (1+\alpha w^{-1})$. In the general
case the result of gluing is the hypersurface
\[
xy-f_{\rho,v}(w_1,\ldots,w_{n-1})t^k=0.
\]

\section{Examples without scattering}

In Section~\ref{sect.SingAffine} we worked in an affine (in the
algebraic geometry sense) neighbourhood of a singular point of the
total space of the degeneration. We now want to look at projective
examples in two dimensions, with one or several singular points. To
this end we work with bounded tropical manifolds and use the
toric procedure from \S\ref{subsect.DegenToric} for the treatment
away from the singularities.

\subsection{Overview of the general procedure}
\label{subsect.overview_general}
To obtain degenerations of complete varieties we need to start from a
bounded tropical manifold $(B,\P)$, where $B$ denotes the underlying
singular integral affine manifold and $\P$ is the decomposition into
integral polyhedra, together with an (in general multi-valued)
integral piecewise affine function $\varphi$. The tropical manifold
determines readily the prospective central fibre $X_0=
\bigcup_{\sigma\in\P_\max} X_\sigma$ from the maximal cells, glued
pairwise torically in codimension one. The most general approach also
allows one to compose the gluing in codimension one with a toric
automorphism, but as this is a rather straightforward generalization
that only complicates the formulas we do not include it here. In the
notation of \cite{affinecomplex} this means taking $s_e=1$ for any
inclusion $e:\omega\to\tau$ of cells $\omega,\tau\in\P$.

To obtain local models near the singular points we need to choose, for
each $\rho\in\P^{[n-1]}$ and vertex $v\in\rho$, a regular function
$f_{\rho,v}$ on an open set in $X_\rho\subseteq X_0$. The open set is
the affine neighbourhood
\[
U_v:=X_0\setminus \bigcup_{\tau\in\P,\,
v\not\in\tau} X_\tau
\]
of the corresponding point $X_v\in X_0$. As explained in
\S\ref{subsect.DegenCodim1}, consistency of the gluing dictates the
change of vertex formula~\eqref{eqn.ChangeOfVertex} that relates
$f_{\rho,v}$ and $f_{\rho,v'}$ for vertices $v,v'\in\rho$. There is
also a compatibility condition for the $f_{\rho,v}$ for fixed $v$ and
all $\rho$ containing a codimension-$2$ cell that makes sure
everything stays toric \'etale locally at general points of the toric
strata. We will explain this condition when we need it later on (see
the discussion after Example~\ref{expl.k=1}).

Next, the piecewise affine function $\varphi$ determines local toric
models for the degeneration, defined by the toric procedure from
\S\ref{subsect.DegenToric}. Explicitly, if $v\in \P$ is a vertex and
$\varphi_v$ is a piecewise linear representative of $\varphi$
introduced in \S\ref{subsect.DegenToric} then the
local model near the zero-dimensional toric stratum $X_v$ is
\[
\Spec\big( \CC[P] \big)\lra \Spec\big(\CC[t]\big),\quad
P=\big\{ (m,h)\in M\times\ZZ\,\big|\, h\ge \varphi_v(m)\big\}.
\]
Here $t=z^{(0,1)}$ makes $\CC[P]$ into a $\CC[t]$-algebra. Following
the discussion in Section~\ref{sect.SingAffine} we now work to some
finite $t$-order $k$, that is modulo $t^{k+1}$, and decompose
according to toric strata. This gives toric local models that depend
on $k$, a toric affine open neighbourhood $U_\omega\subseteq X_0$ of
$\Int(X_\omega)$, a toric stratum $X_\tau$ intersecting $U_\omega$ and
a maximal-dimensional reference cell $\sigma$. In analogy with $U_v$
for $\omega\in\P$ we take
\[
U_\omega:=X_0\setminus \bigcup_{\tau\in\P,\,
\omega\not\subseteq\tau} X_\tau.
\]
Thus we need $\sigma\in\P_\max$ and cells
$\omega,\tau\in\P$ with $\omega\subseteq \tau$. The last condition is
equivalent to $U_\omega\cap X_\tau\neq \emptyset$. The corresponding
ring, localized at all the gluing functions $f_{\rho,v}$ with
$\rho\supseteq \tau$, is denoted $R^k_{\omega\to\tau,\sigma}$, see
\cite{affinecomplex}, \S2.1 for details. Thus $\Spec
R^k_{\omega\to\tau,\sigma}$ is a $k$-th order non-reduced version
(``thickening'') of $(X_\tau\cap U_\omega) \setminus
\bigcup_{\rho\supseteq\tau} V(f_{\rho, v})$, $v\in\omega$ arbitrary.
Note that for $v,v'\in\tau$ the gluing functions $f_{\rho,v}$ and
$f_{\rho,v'}$ differ by a monomial that is invertible on $U_\omega$,
so $V(f_{\rho,v'})\cap U_\omega= V(f_{\rho,v})\cap U_\omega$.

The relation ideal defining $R^k_{\omega\to\tau,\sigma}$ in (the
localization at $\prod_{\rho\supseteq\tau} f_{\rho,v}$ of) $\CC[P]$ is
generated by all monomials $z^{(m,h)}$, $(m,h)\in\Lambda_\tau\oplus
\ZZ$, that have $\tau$-order at least $k+1$. Here one defines the
\emph{$\tau$-order} as the maximum, over all $\sigma'\in\P_\max$
containing $\tau$, of the order of vanishing of $z^{(m,h)}$ on the big
cell of $X_{\sigma'}$, viewed canonically as a subset of $\Spec
\big(\CC[P] \big)$.

Now even in examples with few cells this procedure requires the gluing
of many affine schemes. We can, however, sometimes use the following
shortcut. The integral points of $B$ provide a basis of sections of a
natural very ample line bundle on the central fibre $X_0$, hence
providing a closed embedding $X_0\to \PP^{N-1}$, $N$ the number of
integral points of $B$. In sufficiently simple examples the gluing
morphisms defining the deformation of $X_0$ readily homogenize to
describe the $k$-th order deformation of $X_0$ as a subspace of
$\PP^N_{O_k}$, $O_k=\Spec S_k$, $S_k= \CC[t]/(t^{k+1})$. The toric
nature of this construction makes sure that this procedure gives the
right local models. This is how we compute most of the examples in
this paper.

\subsection{One singular point I: Two triangles}
\label{subsect.TwoTriangles}
\begin{figure}[ht]
\input{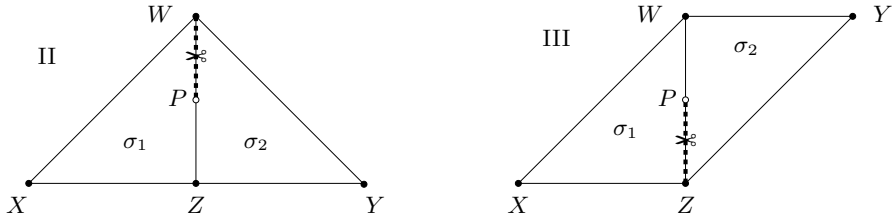}
\caption{The two charts defining the tropical manifold $B_1$.}
\label{fig.B1}
\end{figure}

For our first series of examples we fit the basic example from
\S\ref{subsect.BasicExample} into a projective degeneration. This
amounts to finding a polyhedral decomposition of a bounded
neighbourhood of the singular point $P$ of the focus-focus singularity
in Figure~\ref{fig.focus-focus}. One of the easiest ways to do this is
to take two standard triangles (integrally affine isomorphic to
$\conv\big\{(0,0),(1,0), (0,1)\big\}$), one on each half-plane
$\sigma_\pm$, and with $P$ at the center of the common edge. In
particular, $P$ should not be integral, but rather half-integral.
Thus let us first shift our Charts~I--III by $(0,1/2)$ so that the
singular point $P$ is at $(0,1/2)$. Now in Chart~II take
\[
\sigma_1=\conv\big\{(0,0),(-1,0),(0,1)\big\},\quad
\sigma_2=\conv\big\{(0,0),(1,0),(0,1)\big\}
\]
as maximal cells for the polyhedral decomposition, see
Figure~\ref{fig.B1}. This is one of the two choices making the
boundary locally convex. In fact, we can always bring $\sigma_1$ to
the suggested form by applying a transvection $\big(
\begin{smallmatrix}1&0\\ a&1 \end{smallmatrix} \big)$ on Chart~I for
some $a\in\ZZ$. This descends to an integral affine transformation of
the focus-focus singularity. Then
$\sigma_2=\conv\big\{(0,0),(1,b),(0,1)\big\}$ for some $b\in\ZZ$.
Convexity at $(0,0)$ implies $b\ge 0$, while convexity at $(0,1)$
implies $b\le 1$ (check in Chart~III). Local convexity of $\partial B$
is indispensible to define local toric models for the deformation. 
Both choices are abstractly isomorphic, $b=0$ making the lower
boundary a line, $b=1$ the upper one.
\smallskip

There are $4$~integral points on this tropical manifold $B_1$, the
vertices of $\P$. We interpret these as homogeneous coordinates in a
$\PP^3=\Proj \big( \CC[X,Y,Z,W]\big)$. The central fibre $X_0$ is then
given by the hypersurface $XY=0$, because $X$ and $Y$ are the only
vertices not contained in the same cell. In fact, if $\P$ consists of
standard simplices, $X_0$ can be defined by the Stanley-Reisner ideal
\cite{hochster},\cite{stanley} of the simplicial complex given by
$\P$. This is clear from the definitions.

Now the gluing computation in \S\ref{subsect.BasicExample} readily
homogenizes. In fact, viewing the computation in Chart~II as being
homogenized with respect to Z, the relation \eqref{eqn.BasicExample}
becomes
\[
\frac{X}{Z}\frac{Y}{Z}-t\Big(\alpha+\frac{W}{Z}\Big)=0.
\]
Clearing denominators yields the family of quadric hypersurfaces
\[
XY-tZ(\alpha Z+W)=0.
\]
The analogous computation in the other chart gives the same result. 
Hence we obtain a projective degeneration with general fibre isomorphic
to $\PP^1\times\PP^1$. The total space has another $A_1$-singularity
of a toric nature at $X=Y=Z=0$. This singularity is due to the fact
that $\partial B_1$ is not straight at $(0,1)$. 

It is interesting to compare this deformation to the torically
constructed one for a diagonally subdivided unit square with $\varphi$
changing slope by $1$ along the diagonal.
\begin{figure}[ht]
\input{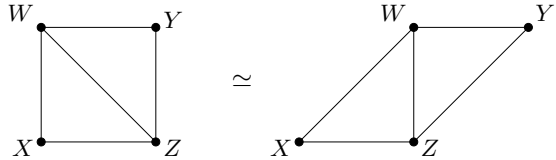}
\caption{Torically constructed degeneration of $\PP^1\times\PP^1$}
\label{fig.DegenP1xP1}
\end{figure}
As shown in Figure~\ref{fig.DegenP1xP1} we again have $4$~homogeneous
generators $X,Y,Z,W$ and the relation
\[
XY-tZW=0.
\]
The change of coordinates $W\mapsto \alpha Z+W$ shows that this family
is isomorphic to the previous one, the only difference being that now
both $A_1$-singularities of the total space are at $0$-dimensional
toric strata of $\PP^3$. Thus in a sense, by introducing the singular
point in the interior we have straightened the boundary of the
momentum polytope at the vertex labelled $Z$ and moved one of the two
$A_1$-singularities in the total space to a non-toric position. The
fact that the introduction of a singular point leads to an isomorphic
family is a rather special phenomenon due to the large symmetry in
this example.

\subsection{One singular point II: Blow-up of $\PP^2$}
\label{subsect:BlownUpP2}
\begin{figure}[ht]
\input{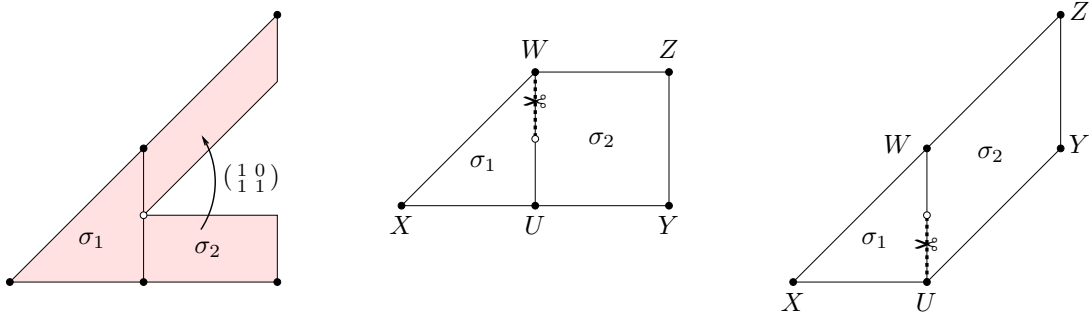}
\caption{Three charts defining the tropical manifold $B_2$.}
\label{fig.B2}
\end{figure}

A similar tropical manifold, which we denote $B_2$, leads to a
degeneration of the blow-up of $\PP^2$ (Figure~\ref{fig.B2}). Note
that the boundary of $B_2$ is straight at both vertices labelled $U$
and $W$. Again we take $\varphi$ to change slope by $1$ along the only
interior edge, with vertices labelled $U$ and $W$. There are
$5$~integral points, so we work in $\PP^4_{\AA^1} =\Proj\big(\CC[t]
[X,Y,Z,U,W] \big)$. 
\medskip
Now $\sigma_2$ is not a standard simplex. Thus we have a relation that
does not receive any corrections from monodromy, namely
\[
WY-UZ=0.
\]
The gluing computations at the vertices labelled $U$ and $W$ homogenize
to the two equations
\begin{eqnarray*}
XY-t(\alpha U+W)U&=&0,\\
XZ-t(\alpha U+W)W&=&0.
\end{eqnarray*}
For $UW\neq0$ these are related by the substitution $Z=WY/U$. The
central fibre $X_0$ has two irreducible components $\PP^2$ and
$\PP^1\times \PP^1$, glued along a $\PP^1$. The total space has only
one singular point, the $A_1$-singularity at $[0,0,0,1,-\alpha] \in
(X_0)_\sing=\PP^1$. The general fibre of this degeneration
$\pi:X\to\AA^1$ contains the $(-1)$-curve $E=V(X,\alpha U+W,\alpha
Y+Z)$ whose contraction yields a $\PP^2$. This is a well-known
example: The $3$~relations form the $2\times2$-minors of the matrix
\[
\left(\begin{matrix} U&W&X\\ Y&Z&t(\alpha U+W)\end{matrix}\right),
\]
which for $t=1$ describe the cubic scroll, see e.g.\ \cite{harris},
Example~7.24.

\begin{figure}[ht]
\input{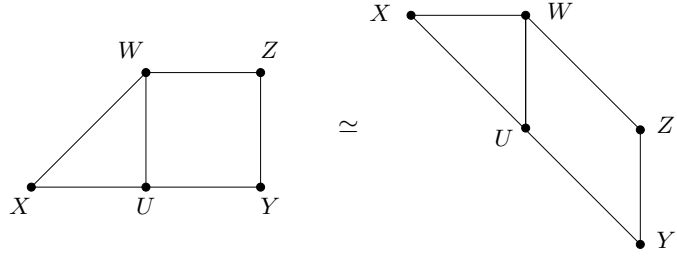}
\caption{Torically constructed degeneration of $\PP^2$ blown up once}
\label{fig.DegenP2BlownUp}
\end{figure}

Again, this example has a toric analogue
(Figure~\ref{fig.DegenP2BlownUp}). The relations are
\[
WY-UZ=0,\quad
XY-tU^2=0,\quad
XZ-tUW=0.
\]
In fact, the substitution $U\mapsto \alpha U+W$, $Y\mapsto \alpha Y+Z$
yields the non-toric ideal above.

\subsection{Propagation}
\label{subsect.Propagation}
Once we start changing the gluing of components somewhere we are
forced to change at other places as well to keep consistency. Thus in a
sense the gluing functions propagate.

\begin{example}
\label{expl.Propagation}
Consider the non-compact example of a tropical manifold with
$4$~maximal cells $\sigma_1,\ldots,\sigma_4$ sown in
Figure~\ref{fig.Propagation}.
\begin{figure}[ht]
\input{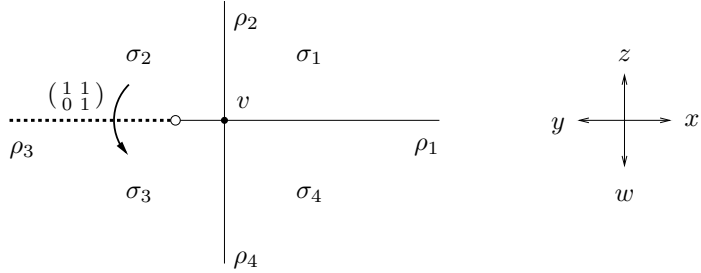}
\caption{A tropical manifold forcing propagation of the gluing
function. The figure on the right shows the tangent vectors belonging
to the generators of the toric model.}
\label{fig.Propagation}
\end{figure}
If we take the piecewise affine function $\varphi$ with
$\varphi(1,0)=1$, $\varphi(0,1)=1$, $\varphi(-1,0)=0$,
$\varphi(0,-1)=0$, the toric local model at $v$ is
\[
\CC[t]\lra \CC[x,y,z,w]/(xy-zw),\quad
t\longmapsto xy=zw.
\]
In terms of integral tangent vectors at $v$ the generators are
\[
x=z^{(1,0,1)},\quad
y=z^{(-1,0,0)},\quad
z=z^{(0,1,1)},\quad
w=z^{(0,-1,0)}.
\]
Now let us try to glue together the $k$-th order neighbourhood of
$X_0$ following the general procedure of \S\ref{subsect.DegenCodim1}
with gluing function $f_{\rho_3,v}=1+y$ and all other
$f_{\rho_i,v}=1$. We have $4$~rings $R^k_{v\to\sigma_i,\sigma_i}$
defining the thickenings of the irreducible components $X_{\sigma_i}$,
$4\cdot 2=8$ rings for the one-dimensional toric strata
$R^k_{v\to\rho_i, \sigma_i}$, $R^k_{v\to\rho_i,\sigma_{i-1}}$ ($i$
taken modulo~$4$), and $4$ rings, identified mutually via an affine
chart at $v$,
\[
R^k_{v\to v,\sigma_i}=\CC[x,y,z,w,t]/
\big( (xy-zw,xy-t)+I_k\big)
\]
for the thickening of the zero-dimensional toric stratum $X_v$. Here
$I_k=(x,z)^{k+1}+(x,w)^{k+1}+ (y,z)^{k+1}+ (y,w)^{k+1}$ is the ideal
generated by monomials that are divisible by $t^{k+1}$ at the generic
point of $X_{\sigma_i}$ for some $i$. For the toric local model we
want to take the inverse limit of all these rings with respect to the
homomorphisms that we have introduced between them. These
are of two kinds. First, for $\tau'\subseteq\tau$ we have $R^k_{v\to
\tau,\sigma_i}\to R^k_{v\to\tau',\sigma_i}$ defining the inclusion of
toric strata with reference cell $\sigma_i$. Second, there are the
change of reference cell isomorphisms $R^k_{v\to\tau,\sigma_i}\to
R^k_{v\to\tau,\sigma_{i\pm 1}}$. This requires compatibility of the
compositions. In the present case this comes down to checking the
following.  Let $\overline{\theta_v^k}$ be the composition
\[
R^k_{v\to v,\sigma_1} \lra
R^k_{v\to v,\sigma_2} \lra
R^k_{v\to v,\sigma_3} \lra
R^k_{v\to v,\sigma_4} \lra
R^k_{v\to v,\sigma_1}
\]
of changing the reference cell $\sigma_i$ counterclockwise around the
origin by crossing the $\rho_i$. The compatibility condition is
$\overline{\theta_v^k}=\id$. But crossing $\rho_1$, $\rho_2$ and
$\rho_4$ yields the identity. Thus $\overline{\theta_v^k}$ equals the
change of reference cell isomorphism from crossing $\rho_3$, which is
\[
x\mapsto x,\quad
y\mapsto y,\quad
z\mapsto (1+y)z,\quad
w\mapsto (1+y)^{-1} w.
\]
Recall that we pick up a negative power of $f_{\rho_3,v}=1+y$ if the
tangent vector points from $\rho_3$ into $\sigma_3$, while monomials
with exponents tangent to $\rho_3$ are left invariant. Note also that
because $y$ is nilpotent now, $1+y$ is invertible, so this is a
well-defined automorphism of $R^k_{v\to v,\sigma_1}$. In any case, we
see $\overline{\theta_v^k}\neq \id$ as soon as $k\ge 2$.

To fix this it is clear how to proceed: Just take $f_{\rho_1,v}=1+y$
as well. Then by the sign conventions, crossing $\rho_1$ gives the
inverse of what we had for $\rho_3$. Then $\overline{\theta_v^k}=\id$ as
automorphism of $R^k_{v\to v,\sigma_i}$ for any $k$ (and any $i$).
\qed
\end{example}

\begin{remark}
There is an important observation to be made here. On
$X_{\rho_3}$ the monomial $y$ is nonzero, that is, of order zero with
respect to $t$. But once we move to a general point of $X_{\rho_1}$ we
have $y=x^{-1}t$ because of the relation $xy=t$. Thus the $t$-order
increases. And this is not an accident, but happens whenever we follow
a monomial $z^{(m,h)}$, $(m,h)\in\Lambda_v\times \ZZ$ in direction
$-m$ (\cite{affinecomplex}, Proposition~2.6). More precisely, if
$\varphi$ differs by $a_\rho \in\Lambda_\rho^\perp$ along a
codimension one cell $\rho$ then a monomial with tangent vector $m$
passing through $\rho$ changes $t$-order by $-\langle
m,a_\rho\rangle$.

This implies that even in examples where the invariant line $\ell$
emanating from a singular point is unbounded, the gluing functions
always converge $t$-adically. Said differently, for the construction
of a $k$-th order deformation of $X_0$, we need to propagate the
contribution of each singular point only through finitely many cells.
\qed
\end{remark}

Of course, there is no reason for the propagation to stay in the
$1$-skeleton of $\P$. The gluing functions then separate from the
$(n-1)$-skeleton of $\P$ and move into the interiors of maximal
cells.

\begin{example}
Let us modify Example~\ref{expl.Propagation} as shown in
Figure~\ref{fig.FirstWallInsertion} with $\varphi(-1,0)=0$,
$\varphi(0,-1)=0$, $\varphi(1,1)=1$. The gluing functions
are $f_{\rho_2,v}=1+x$ and all other $f_{\rho_i,v}=1$.
\begin{figure}[ht]
\input{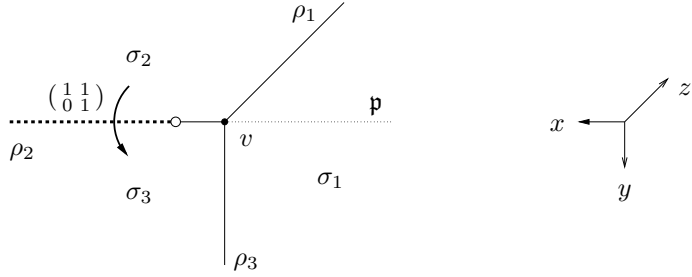}
\caption{Propagation of the gluing function into a maximal cell. The
figure on the right shows the tangent vectors belonging to the
generators of the toric model.}
\label{fig.FirstWallInsertion}
\end{figure}
The toric local model at $v$ is now the normal crossing degeneration
\[
\CC[t]\lra \CC[x,y,z],\quad
t\longmapsto xyz
\]
with $x=z^{(-1,0,0)}$, $y=z^{(0,-1,0)}$, $z=z^{(1,1,1)}$. Following
the gluing isomorphisms around $v$ leads to the automorphism
\[
\overline{\theta_v^k}: x\mapsto x,\quad
y\mapsto (1+x)^{-1} y,\quad
z\mapsto (1+x) z
\]
of the ring
\[
R^k_{v\to v,\sigma_i}=\CC[x,y,z,t]/
\big( xyz-t,x^{k+1},y^{k+1},z^{k+1}\big).
\]
For $k>0$ this is again not the identity, and we can't fix this by
changing any of the gluing functions. But since $1+x$ is invertible on
the whole thickening $\Spec\big(\CC[x,y,z,t]/ ( xyz-t,x^{k+1})
\big)$ of $X_{\sigma_1}$, the inverse of $\overline{\theta_v^k}$ can be
viewed as an automorphism of $X_{\sigma_1}$ that
propagates along the ray $\fop=\RR_{\ge0}\cdot(1,0)$ emanating from
the codimension two cell $v$. We then change the gluing
procedure by taking one copy of the thickening of $X_{\sigma_1}$ for
each connected component of $\sigma_1\setminus\fop$. Passing through
$\fop$ means applying the attached automorphism. The thickenings of
the codimension one strata $X_{\rho_1}$ and $X_{\rho_3}$ are viewed as
closed subschemes of the copies labelled by the connected component of
$\sigma_1\setminus\fop$ containing $\rho_1$ and $\rho_3$,
respectively.

Inserting $\fop$ makes the gluing consistent to all orders. The result
is isomorphic to $\Spec\big( \CC[x,y,z,t]/(xyz-(1+x)t) \big)$.
\qed
\end{example}

Apart from the gluing functions we have now introduced another object,
automorphisms propagating along rays. The higher-dimensional
generalization of rays are \emph{walls} (\cite{affinecomplex},
Definition~2.20): They are one-codimensional polyhedral subsets of
some maximal cell $\sigma$, emanating from a two-codimensional
polyhedral subset $\foq$ into the interior of $\sigma$, and extending
all the way through $\sigma$ along some integral tangent vector
$-m_\fop$, that is, $\fop=\sigma\cap (\foq-\RR_{\ge0}\cdot m_\fop)$.
The attached automorphism is of the form
\[
z^{(m,h)}\longmapsto \big(1+c_\fop z^{(m_\fop,h_\fop)}\big)^{\langle
m,n_{\fop}\rangle}\cdot z^{(m,h)},
\]
where $h_\fop>0$ is the $t$-order of the attached monomial and
$c_\fop\in\CC$. We will discuss how to obtain these walls
systematically in the context of scattering in
Section~\ref{sect.Scattering}.

\subsection{Several singular points: The non-interacting case}
For several singular points $P_1,\ldots,P_l\in B$ the monodromy
invariant direction determines an affine line $\ell_\mu$ emanating
from each $P_\mu$. As long as these coincide or do not intersect, the
construction presented so far works.

\begin{example}
Here is an example with two singular points having parallel invariant
lines $\ell_1$, $\ell_2$.
\begin{figure}[ht]
\input{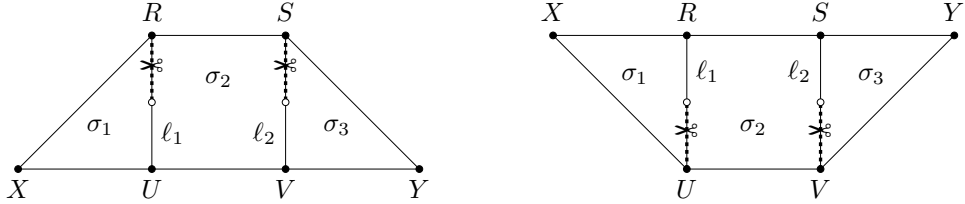}
\caption{Two charts for a tropical manifold with two parallel
invariant lines $\ell_1$, $\ell_2$.}
\label{fig.TwoTrianglesAndSquare}
\end{figure}
We will see that, depending on a choice of parameter, it leads to a
degeneration either of the Hirzebruch surface $\FF_2=\PP
\big(\O_{\PP^1}(-2) \oplus\O_{\PP^1}\big)$ or of $\PP^1\times\PP^1$.

In this case we have two codimension one strata along which we glue.
Taking $\varphi$ to change slope by $1$ along each interior edge
again, the homogenization of the gluing along codimension one cells
gives the polynomials
\[
XV-t(U+ R)U,\quad XS-t(R+U)R,\quad UY-t(V+\lambda S)V,\quad
RY-t(\lambda S+V)S.
\]
Here $\lambda\in\CC^*$ is a parameter that we can not get rid of by
automorphism. Note that $\lambda$ determines the relative position of
the zero loci of the gluing functions when compared in $X_{\sigma_2}
=\PP^1\times\PP^1$. There is one more toric relation $RV-SU=0$ from
$\sigma_2$. The corresponding subscheme of $\PP^5_{\AA^1}=
\Proj_{\CC[t]} \big(\CC[t][X,Y,R,S,U,V] \big)$ has one extra component
$V(R,S,U,V)$. In fact, saturating with respect to any of $R,S,U,V$
gives another relation $XY-t^2(U+R)(V+\lambda S)=0$. This relation can
also be deduced directly from following the line connecting the
vertices $X$ and $Y$ in the chart shown in
Figure~\ref{fig.TwoTrianglesAndSquare} on the right.

These six polynomials are the $2\times 2$ minors of the matrix
\[
\begin{pmatrix}
X&t(V+\lambda S)&U&R\\
t(U+R)&Y&V&S
\end{pmatrix}
\]
For fixed $t\not=0$, this is a scroll in $\PP^5$ whose  ruling
is given by the lines whose equations are linear combinations of the
two rows of the above matrix. The image of this scroll under the
rational map $\kappa:X_t\to \PP^3$ given by $[x,y,u,v,r,s]\mapsto
[x_0,x_1,x_2,x_3]=[u+r,y,v,s]$ has equation $X_0X_1-t(X_2+X_3)
(X_2+\lambda X_3)=0$ in $\PP^3$. For $\lambda\neq 1$ this is a smooth
quadric and $\kappa$ induces an isomorphism $X_t\simeq
\PP^1\times\PP^1$. For $\lambda=1$ the quadric $X_0X_1-t(X_2+X_3)^2=0$
has an $A_1$-singularity, $\kappa$ contracts the
$(-2)$-curve given by $Y=V+S=U+R=X=0$, and $X_t\simeq \FF_2$.
\end{example}

\begin{example}
This example features two singularities with the same invariant line.
\begin{figure}[ht]
\input{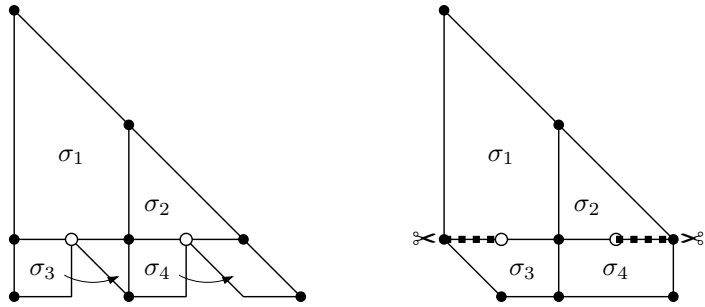}
\caption{A tropical manifold for a degeneration of $\PP^2$ blown
up twice. The invariant lines emanating from the two singular points
agree.}
\end{figure}
The general fibre of the corresponding degeneration is a $\PP^2$ blown
up in two different points. We don't bother to write down the
homogeneous equations.

This example is easy to generalize to any number of
blown up points, which then gives our first genuinely non-toric
examples.
\end{example}

Many more examples can be obtained from \cite{symington}, which
contains a toolkit for the construction of two-dimensional affine
manifolds with focus-focus singularities. The symplectic $4$-manifold
constructed in this reference from such a singular affine manifold is
a symplectic model for the general fibre of our degeneration.

\section{Scattering}
\label{sect.Scattering}
If the invariant lines emanating from singular points of a
two-dimensional affine structure intersect, the gluing construction
becomes inconsistent again.

\subsection{First example of scattering}
\label{subsect.FirstScattering}
\begin{figure}[ht]
\input{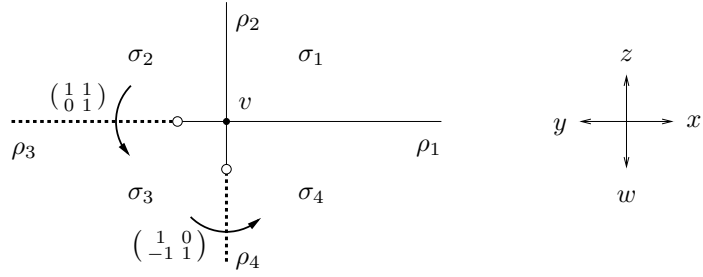}
\caption{A tropical manifold leading to inconsistent gluing. The
figure on the right shows the tangent vectors belonging to the
generators of the toric model.}
\label{fig.FirstScattering}
\end{figure}

We study the modification of Example~\ref{expl.Propagation}
shown in Figure~\ref{fig.FirstScattering}, with gluing functions
\[
f_{\rho_1,v}=1+y=1+x^{-1}t,\quad
f_{\rho_2,v}=1+w=1+z^{-1}t,\quad
f_{\rho_3,v}=1+y,\quad
f_{\rho_4,v}=1+w.
\]
Following the gluing isomorphisms of $R^k_{v\to
v,\sigma_i}$ by going counterclockwise around $v$, with starting cell
$\sigma_1$ say, yields
\begin{eqnarray*}
x&\longmapsto& (1+w)x\\
&\longmapsto& \big(1+(1+y)^{-1}w\big)x\\
&\longmapsto& \big(1+\big(1+(1+w)y\big)^{-1}w\big)(1+w)^{-1}x\\
&\longmapsto&\Big(1+\big(1+(1+(1+y)w)y\big)^{-1}(1+y)w\Big)
\big(1+(1+y)w\big)^{-1}x\\
&=&\Big(1+\big((1+y)+(1+y)wy\big)^{-1}(1+y)w\Big)
\big(1+(1+y)w\big)^{-1}x\\
&=&\big(1+(1+wy)^{-1}w\big) (1+w+wy)^{-1}x\\
&=& \big((1+wy)+w\big)(1+wy)^{-1} (1+w+wy)^{-1}x
\ =\ (1+wy)^{-1}x,
\end{eqnarray*}
and similarly,
\[
y\longmapsto (1+wy)y,\quad
z\longmapsto (1+wy)z,\quad
w\longmapsto (1+wy)^{-1}w.
\]
Thus we can again achieve $\overline{\theta_v^k}=\id$ for all $k$ by
inserting the ray $\fop=\RR_{\ge0}\cdot(1,1)$ with attached function
$1+wy=1+t^2 x^{-1}z^{-1}$!

\subsection{A projective example}
\begin{figure}[ht]
\input{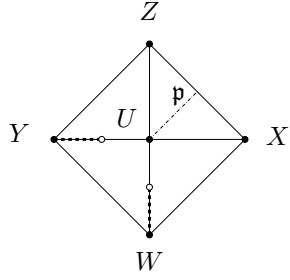}
\caption{Chart for the tropical manifold $B_3$. The dashed lines
denote cuts.}
\label{fig.DiamondTwoSing}
\end{figure}

Of course, this example can also be made projective, say as in
Figure~\ref{fig.DiamondTwoSing}.
The ideal defining the degeneration is generated by the two
polynomials
\[
XY-t(U+W)U,\quad
ZW-t(U+Y)U.
\]
These relations can be most simply derived from gluing $k$-th order
standard thickenings of $X_{\sigma_3}$, $X_{\sigma_4}$ and of
$X_{\sigma_2}$, $X_{\sigma_3}$, respectively, and let $k$ tend to
infinity. This procedure does not seem to depend on $\fop$. In fact,
in this example, the insertion of $\fop$ merely assures that the same
relations are obtained with any other pairs of neighbouring maximal
cells.

As an abstract deformation this is isomorphic to the toric one with
ideal $(XY-tU^2,ZW-tY^2)$ presented in Example~\ref{expl.degenWP1xP1}. 

\subsection{Systematic procedure for the insertion of walls}
\label{subsect.ScatteringProcedure}
The long computation for the insertion of the ray $\fop$ in
Example~\ref{subsect.FirstScattering} looks like quite an accident.
But there is a systematic procedure to insert rays (or walls in higher
dimension) to achieve consistency to any finite order. The number of
rays to be inserted becomes arbitrarily large with increasing $k$,
however, essentially with the only exception the example discussed in
\S\ref{subsect.FirstScattering}. This procedure originates from the
two-dimensional situation of \cite{ks}. Of course, due to the
non-discrete nature of the discriminant locus the higher dimensional
situation is much more involved.

Systematically we proceed as follows. By induction one arrives at a
collection of walls and modifications of the gluing functions, a
\emph{structure} (in its architectural meaning) $\scrS$ as we call it.
A structure consists of walls and so-called \emph{slabs}. A slab
$\fob$ is a polyhedral subset of some $\rho\in\P^{[n-1]}$ together
with a collection $f_{\fob,v}$ of higher order corrections of the
gluing functions $f_{\rho,v}$, one for each vertex $v\in\rho$ that is
contained in a connected component of $\rho\setminus\Delta$
intersecting $\fob$. These fulfill the same change of vertex
formula~\eqref{eqn.ChangeOfVertex} as the $f_{\rho,v}$, but
interpreted in the rings $R^k_{\rho\to\rho,\sigma}$ for any $\sigma\in
\P_\max$ containing $\rho$. We need slabs because a wall hitting a
codimension one cell $\rho$ may lead to different corrections on both
sides of the wall. Technically we refine $\P$ to a polyhedral
decomposition $\P_\scrS$ in such a way that the walls lie in the
$(n-1)$-skeleton of $\P_\scrS$. Then the underlying subsets of the
slabs are exactly the $(n-1)$-cells of $\P_\scrS$ contained in the
$(n-1)$-skeleton of $\P$.

For consistency of the gluing construction it suffices to follow loops
around the codimension two cells of $\P_\scrS$, the \emph{joints} of
the structure (or rather of $\P_\scrS$). This is shown in
\cite{affinecomplex}, Lemma~2.30; the precise definition of
consistency is in Definition~2.28. The fact that the gluing
construction for a structure consistent to order $k$ gives a $k$-th
order toric degeneration occupies \S2.6 in \cite{affinecomplex}.

Thus consistency is really a codimension two feature. For the
scattering computation at a joint $\foj$ we can therefore work over
the ring $S=\CC\lfor t\rfor[\Lambda_\foj]$ of Laurent polynomials with
coefficients in $\CC\lfor t\rfor$ and exponents in $\Lambda_\foj$,
the lattice of integral tangent vectors along $\foj$. In fact, if
$\overline{\theta_\foj}$ is the composition of the gluing morphisms
following a small loop around such a joint, then
$\overline{\theta_\foj}$ is the identity on any monomial $z^{(m,h)}$
with $m$ tangent to $\foj$.

By making the discriminant locus $\Delta$ sufficiently non-rational or
by deforming it a little, we can assume the intersection of $\Delta$
with $\foj$ to be of codimension one in $\foj$. Let
$p\in\foj\setminus\Delta$ be a general point and denote by
$\sigma_\foj\in\P$ the smallest cell containing $\foj$. Note that
$\sigma_\foj$ has dimension $n-2$, $n-1$ or $n$. The computations take
place in an affine chart at a vertex $v\in\sigma_\foj$ in the same
connected component of $\sigma_\foj\setminus \Delta$ as $p$. In an
affine chart at $v$ decompose the lattice of integral tangent vectors
by choosing a complement $\overline \Lambda\simeq\ZZ^2$ to
$\Lambda_\foj \subseteq \Lambda_v$
\[
\Lambda_v=\Lambda_\foj\oplus \overline\Lambda.
\]
While the procedure essentially is always the same, the setup for
the scattering differs somewhat depending on $\dim \sigma_\foj\in
\{n-2,n-1,n\}$. For being explicit let us assume
$\dim\sigma_\foj=n-2$, which is is the most difficult case; it is also
the situation for the initial scattering of codimension one cells of
$\P$ intersecting in a codimension two cell $\foj=\sigma_\foj$ of
$\P$. At a later stage ($k>0$) the inclusion $\foj\subseteq
\sigma_\foj$ may be strict and there may be additional walls running
into $\foj$ from scatterings at other joints.

The scattering procedure is a computation in the rings
$R^k_{\sigma_\foj \to\sigma_\foj,\sigma}$ using induction on $k$.
However, these rings are not practical for computations with a
computer algebra system, and they would require the introduction of
log structures in Example~\ref{expl.k=1}. Log structures are necessary
to make things work globally in the end, but on the level of this
survey it does not seem appropriate to get into these kinds of
technicalities. We therefore present the scattering computation in a
traditional polynomial ring.

For any order $k'\ge 0$ the composition of gluing morphisms around
$\foj$, starting at a reference cell $\sigma\in\P_\max$, defines an
automorphism $\overline{\theta_\foj^{k'}}$ of the ring
$R^{k'}_{\sigma_\foj\to \sigma_\foj,\sigma}$.  This ring is a finite
$S_f$-algebra, where $f= \prod_{\rho\supseteq \sigma_\foj}
f_{\rho,v}$. Let $\overline\varphi$ be the restriction to $\{0\}\times
\overline\Lambda\simeq\ZZ^2$ of a non-negative representative of
$\varphi$. By changing the representative by an affine function we can
assume
\[
\overline\varphi(0,0)= \overline\varphi(-1,0)=0,\quad
0\le \overline\varphi(0,1)\le \overline\varphi(1,1),
\]
which by convexity and integrality of $\overline\varphi$ implies
\[
\overline\varphi(1,1)\ge 1.
\]
Letting $x:=z^{(-1,0,0)}$, $y:=z^{(0,-1,0)}$, $z:=z^{(1,1,1)}$ this exhibits
$R^{k'}_{\sigma_\foj\to \sigma_\foj,\sigma}$ as a quotient of an
$S_f$-subalgebra of
\begin{eqnarray}
\label{R^k}
R^{k}:=S_f[x,y,z]/\big((xyz-t)+(x,y,z)^{k+1}\big),
\end{eqnarray}
for some $k\ge k'$. The smallest possible $k$ is the minimum of the
degrees $a+b+c$ for $x^ay^bz^c$ having $\sigma_\foj$-order (introduced
in~\S\ref{subsect.overview_general}) larger than $k'$. Conversely, any
computation in $R^k$ can be obtained from a computation in
$R^{k'}_{\sigma_\foj\to \sigma_\foj,\sigma}$ for $k'\gg0$. Thus the
result of the scattering procedure is the same with both sorts of
rings.

Initially the only scattering happens at codimension two cells
$\tau\in\P$. The starting data are the gluing functions $f_{\rho,v}$
for all codimension one cells $\rho\supseteq\tau$ where $v\in\tau$ is
some vertex selecting an affine chart for the computation. Consistency
for $k=1$ already requires a compatibility condition for the
restrictions $f_{\rho,v}|_{X_\tau}$ that we can see by studying the
following local example.

\begin{example}\label{expl.k=1}
To consider a local scattering situation for $k=1$ we take
$\foj=\tau=\RR^{n-2} \times\{0\}\subseteq \RR^n$ as an affine subspace
and $\rho_i=  \tau \times \RR_{\ge0} m_i$, $i=1,\ldots,l$, affine
half-hyperplanes emanating from $\tau$ in directions
$m_i\in\ZZ^2=\overline\Lambda$ in the normal lattice $\ZZ^n
/\Lambda_\tau=\overline\Lambda$. (Strictly speaking this is not a
legal situation because $\tau$ does not have a vertex, but this is
irrelevant for the present discussion.) We assume the $\rho_i$
labelled cyclically around $\tau$, that is, $\rho_i$ separates two
maximal cell $\sigma_{i-1}$, $\sigma_i$ containing $\tau$, for any
$i\in \ZZ/l\ZZ$. Let $n_i=(a_i,b_i)$ be a primitive normal vector to
$\Lambda_{\rho_i}$ that evaluates non-negatively on $\sigma_i$.
Then passing from $\sigma_{i-1}$ to $\sigma_i$ gives the automorphism
\[
\overline{\theta^1_i}: z^m\longmapsto f_{\rho_i}^{-\langle
n_i,m\rangle}\cdot z^m
\]
of the ring $R^1=S_f[x,y,z]/(x,y,z)^2$ from \eqref{R^k}. Because in
this local example there is no vertex we drop the reference vertex
from the notation for the gluing functions $f_{\rho_i}$. The loop
around the joint $\tau$ gives the automorphism
$\overline{\theta^1_\tau} =\overline{\theta^1_k} \circ\ldots
\overline{\theta^1_1}$. Because we work modulo $(x,y,z)^2$ the effect
of applying $\overline{\theta^i_k}$ to a monomial is multiplication
with a power of the restriction $f_{\rho_i}|_{X_\tau}\in S_f$.
We obtain
\[
\overline{\theta^1_\tau}(x) = \Big( \prod_i f_{\rho_i}^{-\langle
(a_i,b_i), (-1,0)\rangle}\Big)\cdot x
= \Big( \prod_i \big(f_{\rho_i}|_{X_\tau}\big)^{a_i}\Big)\cdot x,
\]
and analogously,
\[
\overline{\theta^1_\tau}(y) =
\Big( \prod_i \big(f_{\rho_i}|_{X_\tau}\big)^{b_i}\Big)\cdot y,\quad
\overline{\theta^1_\tau}(z) =
\Big( \prod_i \big(f_{\rho_i}|_{X_\tau}\big)^{-a_i-b_i}\Big)\cdot z.
\]
Thus consistency for $k=1$ requires the following multiplicative
condition for the slab functions
\begin{eqnarray}\label{multiplicative condition}
\prod_i \big(f_{\rho_i}|_{X_\tau}\big)^{b_i}=1\quad\text{and}\quad
\prod_i \big(f_{\rho_i}|_{X_\tau}\big)^{a_i}=1.
\end{eqnarray}
\vspace{-5ex}

\qed
\end{example}

Example~\ref{expl.k=1} shows that apart from the change of vertex
formula \eqref{eqn.ChangeOfVertex}, slab functions need to fulfill a
number of algebraic equations to achieve consistency for $k=1$. More
precisely, each $\tau\in\P^{[n-2]}$ yields two multiplicative
equations as in \eqref{multiplicative condition} for the restrictions
$f_{\rho,v}|_{X_\tau}$ of the slab functions for those
$\rho\in\P^{[n-1]}$ that contain $\tau$, see Equation~(1.10) in
\cite{affinecomplex} for the general form. Note that by the change of
vertex formula the relations obtained from different vertices
$v,v'\in\tau$ are equivalent.

For those readers familiar with \cite{logmirror}, what this shows is
that specifying the gluing functions $f_{\rho,v}$ for $k=1$ is
essentially the same as specifying a log structure on $X_0$ (cf.\
\cite{logmirror}, Theorem~3.22). Thus the log structures of
\cite{logmirror} on $X_0$ can be viewed as ``initial conditions'' for
our construction. \medskip

To explain the next steps ($k>1$) let us look at the following typical
scattering situation in dimension~$3$. Let the reference vertex be
$v=0\in\RR^3$ and let the codimension two cell $\tau$ be contained in the
coordinate axis $\RR\cdot (0,0,1)$. We have three codimension one
cells $\rho_1,\rho_2,\rho_3 \supseteq\tau$ with $K_\tau\rho_i/T_\tau$
spanned by $(-1,0,0)$, $(-1,-3,0)$ and $(2,3,0)$, respectively. For
each $i$ there is a unique slab $\fob_i\subseteq \rho_i$ containing
$x$. The piecewise affine function $\varphi$ is uniquely defined by
$\varphi|_\tau=0$ and
\[
\varphi(-1,0,0)=\varphi(-1,-3,0)=0,\quad
\varphi(2,3,0)=3.
\]
Following the maximal cells $\sigma_1, \sigma_2, \sigma_3$ containing
$\tau$ counterclockwise, with $\sigma_1$ containing $\rho_1$ and
$\rho_2$, we have $\varphi(a,b)=0$, $3a-b$ and $b$ for
$(a,b)\in\sigma_1,\sigma_2,\sigma_3$, respectively. In
particular, $\varphi(2,3,0)=3$ is the smallest possible value making
$\varphi$ integral and with $\varphi|_{\sigma_1}=0$. The
exponents of our monomials take values in
\[
P:= \big\{ (m,h)\in\Lambda\times\ZZ\,\big|\,
h\ge\varphi(m) \big\}.
\]
As slab functions take
\begin{eqnarray*}
f_{\fob_1,v}&=& 1+z^{(0,0,1,0)}+ 2z^{(-1,0,0,1)},\\
f_{\fob_2,v}&=& 1+z^{(0,0,1,0)}-z^{(-1,-3,1,0)},\\
f_{\fob_3,v}&=& 1+z^{(0,0,1,0)}+ 5z^{(2,3,0,3)},
\end{eqnarray*}
where the exponents lie in $\Lambda\times\ZZ$, $\Lambda=\Lambda_v$.
Thus in $z^{(a,b,c,h)}$ the first two coordinates $(a,b)$ determine
the direction in $\Lambda/\Lambda_\tau$, the third coordinate $c$ is
for the toric parameter on $X_\tau$ and $h-\varphi(a,b)$ determines
the divisibility by $t$.

\begin{figure}[ht]
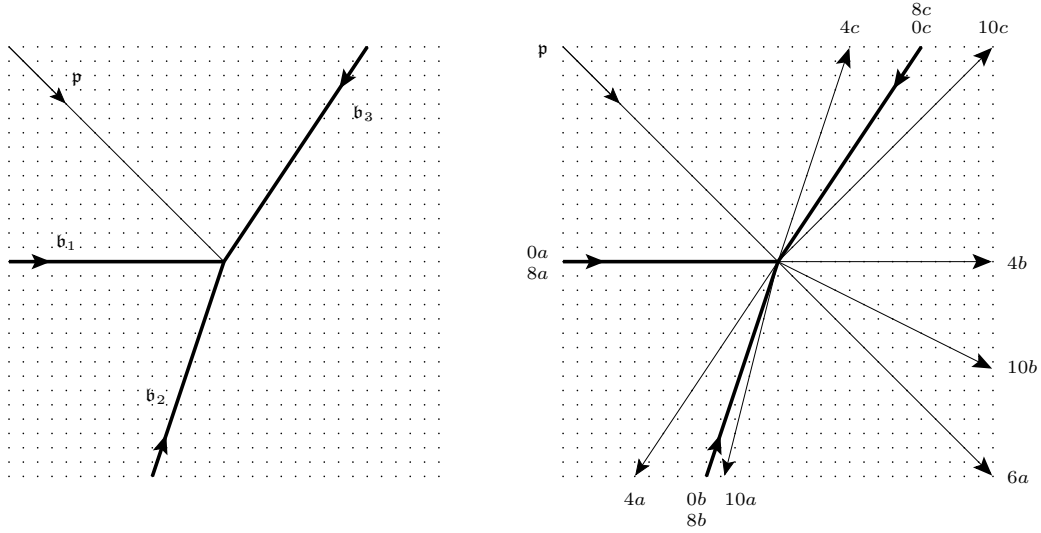
\label{fig.Scattering_Problem}
\input{ScatteringProblem.pstex_t}
\input{ScatteringResult.pstex_t}
\caption{A scattering procedure: Starting data (left) and result
at order $10$. The labels are explained in the text.}
\end{figure}

We also assume there is one wall $\fop$ containing
$\tau$ coming from some earlier step, say with $K_\tau\fop/ T_\tau$
spanned by $(-1,1,0)$ and with associated function
\[
f_\fop=1+7z^{(-1,1,0,2)}.
\]
Our scattering problem can be visualized in the two-dimensional normal
space to $\tau$ as in Figure~\ref{fig.Scattering_Problem} on the left.

For the explicit scattering computation we now go over
to the $S_f$-algebra $R^k$ from \eqref{R^k}. Here $S_f=\CC\lfor t\rfor
[w,w^{-1}]_{1+w}$ with $w=z^{(0,0,1,0)}$, and
\[
x=z^{(-1,0,0,0)},\quad
y=z^{(0,-1,0,0)},\quad
z=z^{(1,1,0,1)}.
\]
Note that $xyz=t$ defines the $t$-algebra structure, and a monomial
$x^a y^b z^c w^d$ represents an element of $P$ if and only if
\[
\varphi(c-a,c-b,d)=\max\big\{0,c-b,3(c-a)-(c-b) \big\} \ge c.
\]
With these notations our input data is
\[
f_{\fob_1,v}= 1+w+ 2tx,\quad
f_{\fob_2,v}= 1+w-wxy^3,\quad
f_{\fob_3,v}= 1+w+ 5xz^3,\quad
f_\fop=1+7tx^2z.
\]
It is now easy to compute the composition $\overline{\theta_v^k}$ by a
computer algebra system. Seeing the various terms come up at higher
order and getting rid of them by insertion of walls is very
instructive. We therefore encourage the reader to do that and verify the
following computations. Note that as we have seen in
\S\ref{subsect.FirstScattering}, the naive insertion procedure leads
to huge expressions quickly, so a computation in $\QQ[x,y,z,w]$ would
be too slow. It is therefore important to expand in a Taylor series to
order $k$ in each of $x$, $y$ and $z$ after each application of an
automorphism associated to a slab or wall. Note also that to compute
$\overline{\theta_v^k}$ as an automorphism of $R^k_{\tau\to\tau,
\sigma_1}$ we have to observe the difference between the $\tau$-order
of a monomial $x^a y^b z^c$ and its degree $a+b+c$. For example, both
$x^2y$ and $x^5y^5z^5=t^5$ have $\tau$-order $5$, but $x^2y$ is
already visible modulo $(x,y,z)^4$ while $x^5y^5z^5$ requires working
modulo $(x,y,z)^{16}$. 

The computation in $R^k$ gives consistency up to $k=3$. At degree
$4$ one has to insert three walls, denoted $4a$, $4b$, $4c$ in
Figure~\ref{fig.Scattering_Problem} on the right. The associated
functions are
\[
1+5(1+w)^2xz^3,\quad
1+2(1+w)^2x^2yz,\quad
1-(1+w)^2wxy^3.
\]
This is just the part $f_{\fob_i}-(1+w)$ of the slab functions that
do not cancel, with an appropriate power of $1+w$ that we will explain
shortly. Similarly, at degree $6$ the incoming wall $\fop$ starts
contributing, forcing its continuation by a wall in direction~$(1,-1)$
with function $1+7(1+w)^5 tx^2z$ ($6a$ in
Figure~\ref{fig.Scattering_Problem} on the right). Then everything
goes well up to degree $8$, where we find interactions of two of such
terms, namely some multiples of $x^3yz^4=tx^2z^3$, $x^3y^4z=t x^2y^3$,
$x^2y^3z^3=t^2 yz$. These monomials have directions $(1,3)$,
$(-2,-3)$, $(1,0)$, so have to scatter away in the directions
$(-1,-3)$, $(2,3)$, $(-1,0)$ occupied by the slabs. Since we do not
want to allow walls to lie above slabs the only possibility is to
change the slab functions. Here is the change that makes things
consistent to degree~$9$:
\begin{eqnarray*}
f'_{\fob_1,v}&=& 1+w+ 2tx-15(1+w)^2t^2wyz ,\\
f'_{\fob_2,v}&=& 1+w-wxy^3+30(1+w)^2tx^2z^3,\\
f_{\fob_3,v}&=& 1+w+ 5xz^3-6(1+w)^2wtx^2y^3.
\end{eqnarray*}
Let us go through degree~$10$ in detail. The composition
modulo $(x,y,z)^{11}$ yields
\begin{eqnarray*}
x&\longmapsto&\Big(1-\frac{700x^4yz^5}{1+w}
-\frac{14x^5y^2z^3}{1+w}-56(1+w)^2wx^4y^4z^2\Big)\cdot x\\
y&\longmapsto&\Big(1+\frac{175x^4yz^5}{1+w}
-\frac{28x^5y^2z^3}{1+w}+56(1+w)^2wx^4y^4z^2\Big)\cdot y\\
z&\longmapsto&\Big(1+\frac{525x^4yz^5}{1+w}
+\frac{42x^5y^2z^3}{1+w} \Big)\cdot z.
\end{eqnarray*}
The occurring monomials are $x^4yz^5=tx^3z^4$, $x^5y^2z^3 =t^2x^3z$,
$x^4y^4z^2=t^2x^2y^2$, which have directions $(1,4)$, $(-2,1)$ and
$(-1,-1)$, respectively. (The fact that the last one points in the
same direction as $z$ explains why it does not show up in the
last line: It belongs to a wall acting trivially on $z$.)
Thus we want to insert walls in directions $(-1,-4)$, $(2,-1)$ and
$(1,1)$ ($10a$, $10b$ and $10c$ in
Figure~\ref{fig.Scattering_Problem}) with functions
\[
1+a tx^3z^4,\quad
1+b t^2x^3z,\quad
1+c t^2x^2y^2,
\]
for some $a,b,c\in\CC[w]$. We now explain how to determine the
coefficients $a,b,c$. Since they come in a product with monomials of
degree $10$, they do not interact with any monomial of non-zero
degree. Hence the influence to the composition can be computed from a
scattering diagram with only the wall in question and the three slabs
$\fob_i$ with functions $f_{\fob_i}=1+w$, $i=1,2,3$. Here is the
computation of the composition for the wall $10a$:
\begin{eqnarray*}
x &&\stackrel{\fob_1}{\longmapsto} \ x
\ \stackrel{\fob_2}{\longmapsto}\  (1+w)^3 x
\ \stackrel{10a}{\longmapsto}\ \big( 1+at x^3 z^4\big)^4\cdot  (1+w)^3 x\\
&&\stackrel{\fob_3}{\longmapsto}\ 
\big( 1+4at (1+w)^{-5}x^3 z^4\big)\cdot x.
\end{eqnarray*}
Here we used $\big( 1+at (1+w)^{-5}x^3 z^4\big)^4 = 1+4at
(1+w)^{-5}x^3 z^4$ modulo $(x,y,z)^{10}$. The second term has to
cancel with $-700 (1+w)^{-1} tx^3z^4$, which leads to $a=175(1+w)^4$.

Said differently, we just have to commute the wall automorphism past
all the automorphisms of slabs in counterclockwise direction up to the
reference cell. For the wall labelled $10a$ there is only one such
slab, $\fob_3$, with normal vector $(-3,2)$, and the monomial is
$x^3z^4=z^{(1,4,0,4)}$. The application of the associated automorphism
yields the power $-\big\langle(-3,2),(1,4)\big\rangle= -5$ of $1+w$.
The coefficient $4=-\big\langle (4,-1), (-1,0) \big\rangle$ is picked
up from passing $x$ past the wall in question. By the form of the
automorphisms, looking at $\overline{\theta_v^k}(y)$ or
$\overline{\theta_v^k}(z)$ gives the same result. This is a
consistency check.

The other coefficients can be computed analogously, giving
\[
1+175(1+w)^4 tx^3z^4,\quad
1+14(1+w)^7 t^2x^3z,\quad
1-56(1+w)^4 t^2wx^2y^2
\]
for the functions associated to the walls $10a$, $10b$ and $10c$,
respectively.
\medskip

Continuing in this fashion we can insert new walls essentially
uniquely to make the gluing consistent to any finite $\tau$-order.

\begin{remark} \label{rem.Scattering}
1)\ There is one miracle here. Namely, while the scattering
computation works only after localization at $1+w$, no negative powers
of $1+w$ ever have to be inserted. This is indispensable because walls
and slabs can only carry regular functions. For general slab functions
and slab configurations this might indeed not be true, see the
discussion in \S4.4.1 of \cite{affinecomplex}. So we need a condition,
and the condition we found to make it work (\cite{affinecomplex},
Definition~1.25,ii) morally says that the deformation problem has a
unique solution locally in codimension~$2$ up to isomorphism. This is
a natural condition. Nevertheless, the proof that no denominators
occur in \S4.4 of \cite{affinecomplex} is the most technical part of
the story. We refer to \cite{affinecomplex} for a more complete
discussion.\\[1ex] 2)\ For $\dim\sigma_\foj\in\{n-1,n\}$ the algorithm
is the same, but there are higher degree monomials with vanishing
$\sigma_\foj$-order. In the $(n-1)$-dimensional case a perturbation
argument as in \S\ref{subsect.perturbation} suffices to show that one
can achieve consistency without changing the slab functions and
without introducing denominators, see \cite{affinecomplex} \S4.3. In
the $n$-dimensional case this is not an issue because there are no
slabs and hence there are no denominators and all monomials have
positive order anyway. This last case is essentially the situation
that was treated by Kontsevich and Soibelman in \cite{ks}. In essence,
our approach trades the difficult gradient flow arguments in \cite{ks}
for algebraic arguments for the $\dim\sigma_{\foj}=n-1$ case, while
the $\dim\sigma_\foj=n-2$ case is only non-trivial for $n\ge 3$ and is
the most difficult issue.
\end{remark}

\subsection{The perturbation idea}
\label{subsect.perturbation}
Sometimes a simple geometric perturbation of the starting data allows
one to simplify a scattering situation. In \cite{affinecomplex}, \S4.2
we formalize this in the notion of \emph{infinitesimal scattering
diagrams}. Here we just want to illustrate the concept by one simple
example.

\begin{example}
Consider the two-dimensional example of a local scattering situation
with $xyz=t$ shown in Figure~\ref{fig.PerturbationData}.
\begin{figure}[ht]
\input{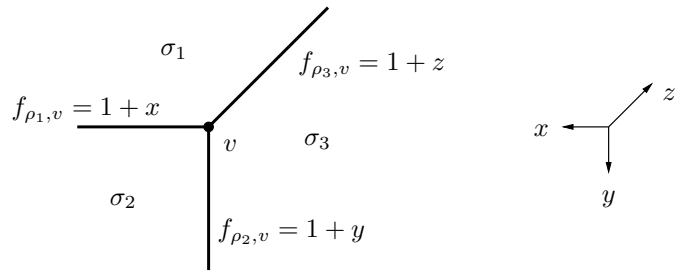}
\caption{Scattering of three non-trivial slabs.}
\label{fig.PerturbationData}
\end{figure}
By formally
perturbing the slabs we arrive at Figure~\ref{fig.Perturbation} on the
left. The functions are the slab functions on the various polyhedral
pieces. While this is not a legal polyhedral decomposition, we can
still insert walls inductively to make the diagram consistent to any
finite order. The result is shown in Figure~\ref{fig.Perturbation} on
the right. Note that this diagram is entirely determined by the most
simple scattering situation of \S\ref{subsect.FirstScattering} via
affine isomorphisms.
\begin{figure}[ht]
\input{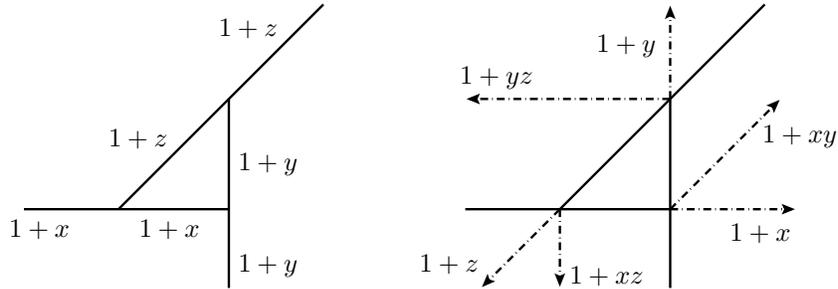}
\caption{Perturbed scattering diagram.}
\label{fig.Perturbation}
\end{figure}

Now the point of this discussion is that we can collapse this diagram
again by superimposing walls and slabs with the same direction if
\begin{figure}[ht]
\input{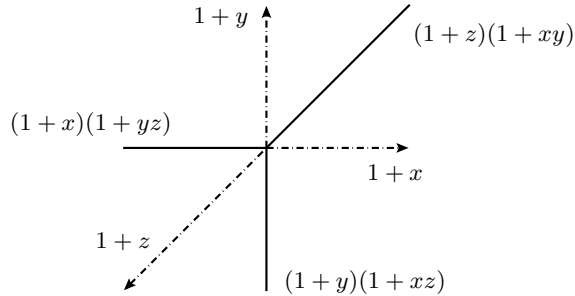}
\caption{Result of collapsing Figure~\ref{fig.Perturbation}, right.}
\label{fig.PerturbationResult}
\end{figure}
necessary. For the present situation the result is shown in
Figure~\ref{fig.PerturbationResult}. This configuration of slabs and
walls is consistent to all orders.
\qed
\end{example}

We can also make this example projective. The smallest
polarization produces a degeneration of cubic
surfaces as follows.

\begin{example}
Consider the tropical manifold depicted in
Figure~\ref{fig.DegenCubic}, with three unit triangles as maximal
cells and one focus-focus singularity on each of the three interior
edges.
\begin{figure}[ht]
\input{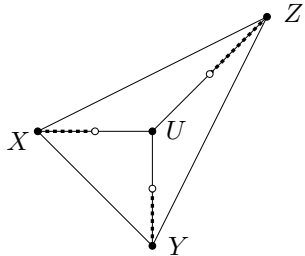}
\caption{A tropical manifold for a degeneration of cubic surfaces.}
\label{fig.DegenCubic}
\end{figure}
The piecewise affine function is defined by the values $0,0,0,1$
on the $4$~vertices $U,X,Y,Z$, respectively. The vertices also define
homogeneous coordinates written by the same symbols. This exhibits
$X_0$ as the union of hyperplanes $V(XYZ)\subseteq \PP^3$.
Dehomogenizing with respect to $U$ gives the affine coordinates
$x_0=X/U$, $y_0=Y/U$, $z_0=Z/U$, and in these coordinates we define
the slab functions by
\[
f_{\rho_1,U}=1+x_0,\quad
f_{\rho_2,U}=1+y_0,\quad
f_{\rho_3,U}=1+z_0.
\]
Let us compute an affine chart for the deformation. To this end we
look at Figure~\ref{fig.PerturbationResult} and perform a computation
similar to the basic gluing computation
in~\S\ref{subsect.BasicExample}. As reference chamber choose the upper
left quadrant in Figure~\ref{fig.PerturbationResult}. In this chamber
$x=X/U$ is the uncorrected toric coordinate $x_0$, but $y=Y/U$ and
$z=Z/U$ receive corrections from crossing walls and slabs. For $z$ we
need to cross the wall with function $1+y_0$, while for $y$ there is the
wall with function $1+z_0$ and the slab with function
$(1+x_0)(1+y_0z_0)$. The
complete list of relations therefore is
\begin{eqnarray*}
x_0y_0z_0&=&t\\
x&=&x_0\\
z&=&(1+y_0)z_0\\
y&=&(1+x_0+z_0+y_0z_0+x_0y_0z_0)y_0.
\end{eqnarray*}
Eliminating $x_0,y_0,z_0$ exhibits an affine chart for the deformation
as the hypersurface
\[
xyz=t\big( (1+t)+x+y+z\big).
\]
in $\AA^4$ with coordinates $x,y,z,t$. Homogenizing yields the
degeneration of singular cubic surfaces
\[
XYZ=t\big( (1+t)U^3+(X+Y+Z)U^2 \big).
\]
The general fibre has $3$~ordinary double points at $[1,0,0,0]$,
$[0,1,0,0]$ and $[0,0,1,0]$; these come from the vertices $X,Y,Z$ on the
boundary of $B$.
\qed
\end{example}

\section{Three-dimensional examples}
So far we have essentially considered two-dimensional examples, the
only exception being the three-dimensional sample scattering
computation in \S\ref{subsect.ScatteringProcedure}. We already
observed in this example one complication in higher dimensions, the
potential presence of poles in the scattering procedure, see
Remark~\ref{rem.Scattering},(1) for a discussion how this is handled.
We also observed in this example that the scattering procedure
generally requires higher order corrections to the slab functions, and
the slab functions themselves propagate. In two dimensions ad hoc
solutions can be used to do this. In higher dimensions this
propagation is more complicated and we need a homological argument
(\cite{affinecomplex}, \S3.5). Another, more fundamental difference in
higher dimensions is the fact that the codimension two intersection
loci of walls and slabs, the joints, are higher dimensional. Thus in
the automorphism $\overline{\theta_\foj}$ associated to a loop around
a joint $\foj$ there may be monomials tangent to $\foj$. These can not
be removed by inserting walls.

It nevertheless turns out that once we have run the homological
argument, this only happens for joints $\foj$ contained in a
codimension two cell $\sigma_\foj$ of $\P$, and the remaining terms
are either undirectional (pure $t$-powers, $t^l$) or of the form $t^k
z^{m^\rho_{vv'}}/ f_{\rho,v}$ (\cite{affinecomplex},
Proposition~3.23). Here $m^\rho_{vv'}$ is defined by affine monodromy
as in \eqref{eqn.ParallelTransport}. In other words, these terms arise
from pure $t$-powers on some other affine chart. The final step in our
algorithm uses a normalization procedure to get rid of these terms
also  (\cite{affinecomplex}, \S3.6). The normalization essentially
says that the logarithms of the slab functions do not contain pure
$t$-powers. Remarkably, this step not only makes sure that no
obstructions arise in the deformation process, but also makes $t$ a
canonical parameter in the sense of mirror symmetry, traditionally
denoted~$q$.

Examples featuring all phenomena in this process are very complicated
to run through explicitly. We therefore content ourselves with one
non-compact example treating a degeneration of the total space of
$K_{\PP^2}$, the so-called ``local $\PP^2$'' from the mirror symmetry
literature \cite{local mirror}, and its mirror.

\begin{example}\label{expl.LocalP2}
Consider the tropical manifold $B$ shown in Figure~\ref{fig.LocalP2}.
\begin{figure}[ht]
\input{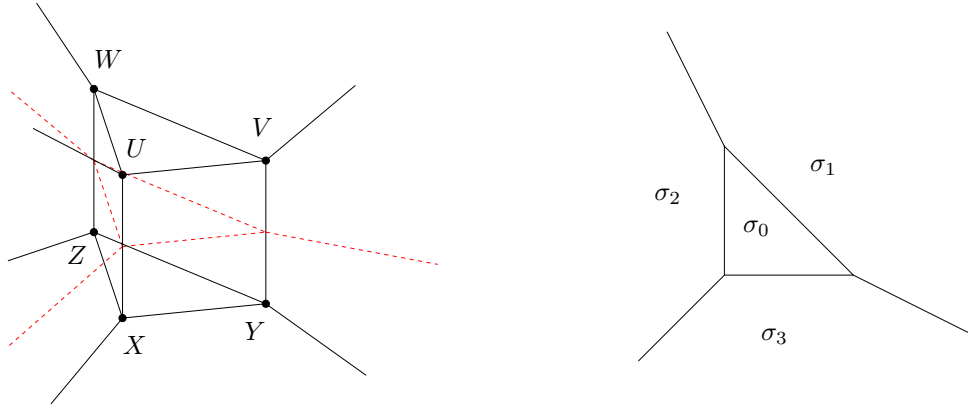}
\caption{The tropical manifold for a degeneration of $K_{\PP^2}$
(left) and its projection to the plane (right).}
\label{fig.LocalP2}
\end{figure}
There is a total of six maximal cells $\sigma_0,\ldots,\sigma_5$ with
only $\sigma_0$ bounded, a prism. The maximal cells adjacent to the
sides of the prism are $\sigma_1,\sigma_2,\sigma_3$ as shown in the
figure. The remaining $\sigma_4,\sigma_5$ are adjacent to the bottom
and top of $\sigma_0$, respectively. The affine structure is then
completely determined by
\[
\sigma_0= \conv\big\{(0,0,0), (1,0,0),(0,1,0),
(0,0,1),(1,0,1),(0,1,1)\big\},
\]
by the directions $(2,-1,\pm1)$, $(-1,2,\pm1)$, $(-1,-1,\pm1)$ of the
unbounded edges in a chart on the upper (plus sign) and lower (minus
sign) half, and by requiring the monodromy around the discriminant
locus (the dashed lines in the figure) to be primitive and positive.
Note that $(0,0,1)$ is an invariant tangent vector defining an affine
projection as suggested in Figure~\ref{fig.LocalP2} on the right.
As piecewise affine function $\varphi$ we take a minimal one changing
slope by $1$ along each codimension~one cell and vanishing on
$\sigma_0$. In particular, on the six unbounded edges $\varphi$
takes the value $1$ at the first integral point different from a
vertex.

The nontrivial slabs ($f_\fob\neq 1$) are the six vertical cells of
codimension one. If we denote by $s$ the monomial with direction
$(0,0,-1)$ and with $\sigma_0$-order $0$ and by $s'$ its inverse, the
gluing functions are
\[
f_\fob=1+s\quad\text{(upper half)},\qquad
f_\fob=1+s'\quad\text{(lower half)}.
\]
Note that by \eqref{multiplicative condition}
and the change of vertex formula \eqref{eqn.ChangeOfVertex}
any one slab function determines all the others.

This structure is already consistent along the three vertical edges of
the prism. To get consistency everywhere it only remains to insert six
vertical walls in $\sigma_4$ and $\sigma_5$ each. These extend the six
vertical codimension one cells to infinity. There is no scattering
because the monomials carried by the walls and slabs all point in the
same direction $(0,0,1)$, the invariant tangent vector, and hence the
corresponding automorphisms mutually commute.

To study this example it is not advisable to write down the
homogeneous coordinate ring. In fact, because this is a non-compact
example it does not suffice to take the six vertices $X,Y,Z,U,V,W$ as
generators of a homogeneous coordinate ring. Rather we need a number
of generators of degree $0$ defined by tangent vectors in the
unbounded directions. Due to the non-simplicial nature of the
polyhedra this leads to a long list of generators and an even longer
list of relations.

Instead we use the construction via gluing affine patches from
\cite{affinecomplex}. There are six vertices, and correspondingly we
have a cover of the degeneration by six affine open sets.
Figure~\ref{fig.LocalP2Coord}
\begin{figure}[ht]
\input{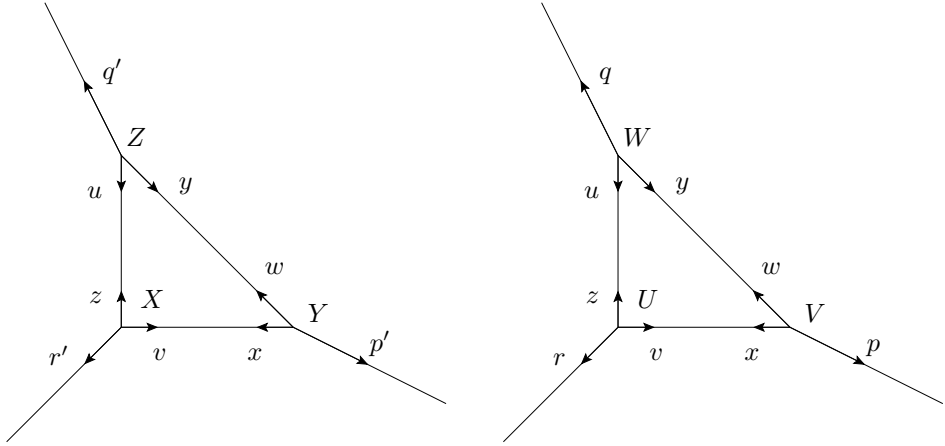}
\caption{The generators of the affine patches (left: bottom, right:
top).}
\label{fig.LocalP2Coord}
\end{figure}
shows our choices of generators, where
the arrows should be thought of as tangent to edges. Thus in a chart
at $U$ the generating monomials are
\[
v=z^{(1,0,0,0)},\ z=z^{(0,1,0,0)},\ r=z^{(-1,-1,1,1)},\ 
s=z^{(0,0,-1,0)},\ t=z^{(0,0,0,1)},
\]
with single relation
\[
rvzs=t,
\]
a semi-stable (normal crossings) degeneration. Because of the symmetry
of the example the situation is analogous in the five other charts.

Now we have to adjust these local models by the slab functions. A
computation analogous to \S\ref{subsect.BasicExample} gives
\begin{align*}
rvzs&=(1+s)t \quad\text{(at $U$)},&
pwxs&=(1+s)t \quad\text{(at $V$)},&
quys&=(1+s)t \quad\text{(at $W$)},\\
r'vzs'&=(1+s')t \quad\text{(at $X$)},&
p'wxs'&=(1+s')t \quad\text{(at $Y$)},&
q'uys'&=(1+s')t \quad\text{(at $Z$)}.
\end{align*}
Here the variables take reference to any maximal cell containing them.
Thus $x,y,z,u,v,w$ can all be thought of as monomials on $\sigma_0$,
while for example $r$ is a monomial on $\sigma_2$ or $\sigma_3$. The
patching between these charts is given by affine geometry in the
relevant maximal cells. In particular, variables with the same name
are all identified, and we have relations such as
\[
x=v^{-1},\quad
u=x w^{-1},\quad
s'=s^{-1},\quad
q=pw^3\quad\textrm{etc.}
\]
As a consistency check one can verify compatibility of the relations
with the gluing. For example, going from the chart at $U$ to the chart
at $V$ means the substitution
\[
v\mapsto x^{-1},\quad z\mapsto x^{-1}w,\quad
r\mapsto x^3p,\quad s\mapsto s.
\]
This maps the relation $rvzs=(1+s)t$ to $pwxs=(1+s)t$, as expected.
Similarly, to go from $U$ to $X$ means substituting
\[
v\mapsto v,\quad z\mapsto z,\quad r\mapsto r's',\quad s\mapsto
(s')^{-1}
\]
into $rvzs=(1+s)t$, leading to $r'vzs'=(1+s')t$.

At this point we have written down a degeneration $\pi:X\to
\AA^1=\Spec \CC[t]$ with $X$ covered by six affine open sets. We claim
that a general fibre $X_t$ is an open subset of the total
space $K_{\PP^2}$ of the canonical bundle of $\PP^2$. To this end
fix $t\in\CC\setminus\{0\}$ and define a projection
\[
\kappa: X_t\longrightarrow \PP^2
\]
by viewing the triples $X,Y,Z$ or $U,V,W$ as homogeneous coordinates
on $\PP^2$. Thus set-theoretically the restriction of $\kappa$ to the
chart at $U$ is
\[
(r,v,z,s)\longmapsto [1,v,z].
\]
It is straightforward to check compatibility with the patching. For
example, in the intersection with the chart at $V$ we find
\[
\kappa(p,w,x,s)= [x,1,w] = [v^{-1},1,v^{-1}z]=[1,v,z].
\]
Analogous computations show compatibility on the ring level.

The fibre of $\kappa$ over a closed point of $\PP^2$, say $[1,v,z]$,
is the hypersurface $rvzs-(1+s)t=0$ in $\AA^2=\Spec\big(
\CC[r,s]\big)$. Note that $X_t$ is disjoint from $s=0$ or from $s'=0$
since $t\neq0$, so it suffices to work in one chart only. If $vz\neq0$
this is a hyperbola, hence isomorphic to $\AA^1\setminus\{0\}$. On the
other hand, if $vz=0$ we must have $s=-1$ and $r$ has no
restrictions, so this is an $\AA^1$. The global meaning of this comes
by observing that the fibre coordinates $p,q,r$ transform dually to
the sections $dx\wedge dw$,
$dy\wedge du$ and $dz\wedge dv$ of $K_\PP^2$. For example, since
$y=z^{-1}v$ and $u=z^{-1}$,
\[
qdy\wedge du=q\big(-z^{-2}vdz +z^{-1}dv)\wedge \big(-z^{-2}dz\big)
= qz^{-3}dz\wedge dv =rdz\wedge dv.
\]
Above we computed the fibre of $\kappa$ over $[1,v,z]$ to be given by
$(rvz-t)s=t$. For given $r,v,z$ this has a solution $s$ as long as
$rvz-t\neq0$. Thus $rvz-t=0$ describes the hypersurface locally that
is being removed from $K_{\PP^2}$ to obtain $X_t$. Globally we are
removing the graph of a rational section of $K_{\PP^2}$ with poles
along the toric divisor $V(XY\!Z)$.
\qed
\end{example}

\begin{example}\label{expl.MirrorLocalP2}
An example with somewhat complementary features to
Example~\ref{expl.LocalP2} is provided by the mirror. Following the
general recipe of \cite{logmirror} the tropical manifold together with
the (multi-valued) function $\varphi$ is obtained by a discrete
version of the Legendre transform (\cite{logmirror},
Construction~1.15). By this construction the polyhedral decompositions
of a tropical manifold and its mirror are combinatorially dual to each
other. As this is not primarily a paper about mirror symmetry we do
not explain this construction, but only state the result.

As shown in Figure~\ref{fig.MirrorLocalP2}
\begin{figure}[ht]
\input{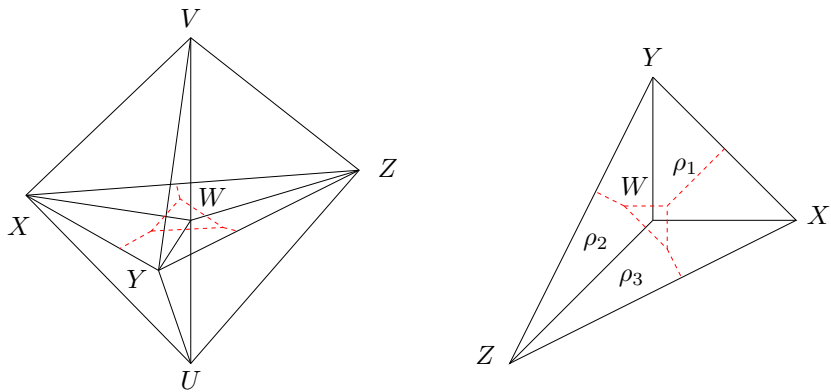}
\caption{The tropical manifold for the mirror of
Example~\ref{expl.LocalP2} (left) and the central triangle containing
the discriminant locus (right).}
\label{fig.MirrorLocalP2}
\end{figure}
the tropical manifold is a double tetrahedron glued from six standard
simplices. We have five exterior vertices labelled $X,Y,Z,U,V$ and one
interior vertex $W$. In an affine chart at $W$ the six emanating edges
point in the directions of the coordinate axes. The affine structure
is then completely determined by requiring the monodromy around the
edges of the discriminant locus (dashed in the figure) to be primitive
and positive. The function $\varphi$ is single-valued and can be taken
to take value $0$ on $X,Y,U, W$ and value $1$ on $Z$ and $V$. This again
has the property to change slope by $1$ along each cell of codimension
one.

There are three non-trivial slabs, the three horizontal triangles
containing the discriminant locus. Up to automorphisms there is only
one set of gluing functions possible at order $0$, namely, in affine
coordinates $x=X/W$, $y=Y/W$, $z=Z/W$ at $W$:
\[
f_{\rho_1,W}=1+x+y,\quad
f_{\rho_2,W}=1+y+z,\quad
f_{\rho_3,W}=1+z+x.
\]
The expressions at the other vertices follow from this by the change
of vertex formula~\eqref{eqn.ChangeOfVertex}.

To make this structure consistent to all orders only requires
propagating the slab functions to the neighbouring slabs. This leads
to
\[
f_{\fob,W}=1+x+y+z,
\]
for any of the three slabs $\fob=\rho_i$.

The approach by homogeneous coordinates works well again in this case.
We get the toric relation
\[
XYZ=tW^3,
\]
and the homogenization of the gluing relation $uv=(1+x+y+z)t$ gives
\[
UV=t^2(X+Y+Z+W)W.
\]
As one checks in local coordinates this is a complete set of
relations. For $t\neq 0$ the projective variety $X_t$ is a conic
bundle over $\Proj\big(\CC[X,Y,Z,W]/(XYZ-tW^3)\big)$ with singular
fibers over $(X+Y+Z+W)W=0$. This base space of the conic bundle is the
quotient of $\PP^2$ by the $\ZZ/3$-action
\[
\xi\cdot[x_0,x_1,x_2]= [\xi x_0,\xi^2 x_1,x_2],
\]
for $\xi$ a primitive third root of unity.

The suggestion in the literature for the mirror to $K_{\PP^2}$ is to
take the family of non-complete Calabi-Yau varieties defined by
\[
uv=1+x+y+tx^{-1}y^{-1}
\] 
in $\CC^2\times(\CC^*)^2$ \cite{local mirror}. Here $u,v$ are the
coordinates on $\CC^2$. This is exactly the open subset
of our family fibering over the big cell of the weighted projective
space. See \cite{Gr1}, \S4 for a discussion how this fits with the SYZ
picture of local mirror symmetry.

As written this family does not come correctly parametrized for the
purpose of mirror symmetry. Rather, a period integral defines a new
parameter $q$ related to $t$ by the so-called \emph{mirror map}. It is
one striking feature of our approach that this mirror map comes up
naturally via the normalization condition. The present example is too
local to illustrate the need for doing this, but as mentioned at the
beginning of this section, generally our algorithm requires the
logarithm of the slab functions to not contain any pure $t$-powers,
see \cite{affinecomplex}, \S3.6. In the present example this means
adding a power series $g=\sum_{l\ge 0} a_l t^l$ to $f_{\fob,W}$ with
the property that
\[
\log\big( f_{\fob,W}+g\big)
=\sum_{k\ge1} \frac{(-1)^{k+1}}{k}\big(x+y+z+g(xyz)\big)^k
\in\CC\lfor x,y,z\rfor
\]
does not contain any monomials $(xyz)^l=t^l$. This condition
determines the coefficients $a_k$ inductively:
\[
g(t)=-2t+5t^2-32t^3+286t^4-3038t^5+\ldots
\]
It follows from the period computations in \cite{GbZa} that the
modified family
\[
XYZ=tW^3,\quad UV=t^2\big(X+Y+Z+(1+g(t))W\big)W,
\]
is then indeed written in canonical coordinates, that is, the mirror
map becomes trivial.
\qed
\end{example}



\begin{thebibliography}{cccccccc}
\bibitem[ChKlYa]{local mirror} T.M.~Chiang, A.~Klemm,
	S.-T.~Yau, E.~Zaslow:
	\emph{Local mirror symmetry: calculations and interpretations},
	Adv.\ Theor.\ Math.\ Phys.~\textbf{3} (1999), 495--565.
\bibitem[GbZa]{GbZa} T.~Graber, E.~Zaslow: 
	\emph{Open-string Gromov-Witten invariants: calculations and 
	a mirror ``theorem''.}
	in ``Orbifolds in mathematics and physics (Madison, WI, 2001)'',  
	107--121, Contemp. Math., {\bf 310}, 
	Amer.\ Math.\ Soc.~2002.
\bibitem[Gr1]{Gr1} M.~Gross:
	\emph{Examples of special Lagrangian fibrations},
	in: \textsl{Symplectic geometry and mirror symmetry}
	(Seoul, 2000),  81--109, World Sci.\ Publ.~2001
\bibitem[Gr2]{GBB} M.~Gross:
	\emph{Toric Degenerations and Batyrev-Borisov Duality},
	Math.\ Ann.\ \textbf{333} (2005), 645--688.
\bibitem[GrSi1]{logmirror} M.\ Gross, B.\ Siebert:
	\emph{Mirror symmetry via logarithmic degeneration data I},
	J.\ Differential Geom.~\textbf{72}
	(2006), 169--338.
\bibitem[GrSi2]{affinecomplex} M.\ Gross, B.\ Siebert:
	\emph{From real affine to complex geometry},
	preprint \texttt{arXiv:math/0703822}, 128~pp.
\bibitem[Ha]{harris} J.\ Harris:
	\emph{Algebraic geometry},
	Springer 1992.
\bibitem[Ho]{hochster} M.\ Hochster:
	\emph{Cohen-Macaulay rings, combinatorics and simplicial
	complexes},
	in: \textsl{Ring theory II}, B.R.~McDonald, R.A.~Morris (eds.),
	Lecture Notes in Pure and Appl.\ Math.~26, M.~Dekker~1977.
\bibitem[KoSo]{ks} M.\ Kontsevich, Y.\ Soibelman:
	\emph{Affine structures and non-Archimedean analytic spaces}, in:
	\textsl{The unity of mathematics} (P.~Etingof, V.~Retakh,
	I.M.~Singer, eds.),  321--385, Progr.\ Math.~244,
	Birkh\"auser~2006.
\bibitem[St]{stanley} R.~Stanley:
	\emph{Combinatorics and commutative algebra, Second ed.},
	Birkh\"auser 1996.
\bibitem[Sy]{symington} M.\ Symington:
	\emph{Four dimensions from two in symplectic topology},
	in: \textsl{Topology and geometry of manifolds (Athens, GA, 2001)},
	153--208, Proc.\ Sympos.\ Pure Math.~71, Amer.\ Math.\ Soc.~2003.
\bibitem[Wi]{williamson} J.~Williamson:
	\emph{On the algebraic problem concerning the normal forms of
	linear dynamical systems},
	Amer.\ J.\ Math.~\textbf{58} (1936), 141--163.
\end{thebibliography}
\end{document}